\documentclass[10pt]{article}

\usepackage[utf8]{inputenc}

\usepackage[margin=1in]{geometry}

\usepackage{hyperref}
\usepackage{optidef}
\usepackage{amsmath,amssymb,amsfonts,amsthm,enumerate}
\usepackage{eucal} 
\usepackage{xcolor}
\usepackage{youngtab}
\usepackage{young}
\usepackage{lscape}
\usepackage{tikz}
\usepackage{environ}
\usepackage{caption}
\usepackage{subcaption}
\usepackage{algorithm}
\usetikzlibrary{positioning, backgrounds}
\usepackage[noend]{algpseudocode}
\usepackage{graphicx} 
\usepackage{amsmath,amsthm,amssymb,amscd}
\usepackage[export]{adjustbox}

\DeclareMathOperator{\mat}{Mat}
\newcommand{\densalpha}{\mathcal{X}}
\newcommand{\densalphao}{\mathcal{X}^{\circ}}

\DeclareMathOperator{\affhull}{aff}

\newcommand{\betti}{\beta}
\newcommand{\quantile}{Q}

\newcommand{\grad}{\nabla}

\newcommand{\levset}{\mathcal{L}}

\DeclareMathOperator{\isingint}{int}
\DeclareMathOperator{\isingcirc}{circ}
\DeclareMathOperator{\isingflares}{flares}

\newcommand{\R}{\mathbb{R}}

\newcommand{\dataset}{\mathcal{D}}

\newcommand{\shape}{\mathcal{S}}

\DeclareMathOperator{\nerve}{Nrv}

\DeclareMathOperator{\spn}{span}

\DeclareMathOperator{\chull}{conv}

\newtheorem*{thma}{Theorem A}

\newtheorem{thm}{Theorem}

\newtheorem{lemma}{Lemma}

\theoremstyle{definition}
\newtheorem{defn}{Definition}
\newtheorem{remark}{Remark}

\theoremstyle{definition}

\newtheorem{alg}{Algorithm}
\newtheorem{prob}{Problem}

\counterwithin{figure}{subsection}

\counterwithin{equation}{section}

\counterwithin{alg}{section}

\title{Alpha shapes in kernel density estimation}

\author{Erik Carlsson and John Carlsson}
\begin{document}
\maketitle
\begin{abstract}
For every Gaussian kernel density estimator
$f(x)=\sum_i a_i \exp(-\lVert x-x_i\rVert^2/2h^2)$
associated to a point cloud 
$\dataset=\{x_1,...,x_N\}\subset \mathbb{R}^d$,
we define a nested family of closed subspaces 
$\shape(a)\subset\mathbb{R}^d$,
which we interpret as a continuous version of an \emph{alpha shape}.
Using arguments based on Fenchel duality,
we prove that $\shape(a)$ is homotopy equivalent to the superlevel 
set $\levset(a)=f^{-1}[e^{-a},\infty)$, and that $\levset(a)$
can be realized as the union of a certain power-shifted
covering by balls with centers in $\shape(a)$. 
By extracting finite alpha complexes with vertices in 
$\shape(a)$, we obtain refined geometric models of noisy point clouds, as well as density-filtered persistent homology calculations.
In order to compute alpha complexes in higher dimension, we used a recent algorithm due to the present authors based on the duality principle \cite{carlsson2023alpha}.
\end{abstract}

\section{Introduction}

Let $f:\mathbb{R}^d\rightarrow \mathbb{R}_+$ be a sum of Gaussian kernels
with uniform covariance matrices, which after 
a linear change of coordinates takes the form
\begin{equation}
\label{eq:defkde}
    f(x)=\sum_{i=1}^N a_i \exp\left(-\lVert x-x_i\rVert^2/2h^2\right),
\end{equation}
where $a_i>0$, $h>0$ is called the bandwidth or scale parameter,
and $\dataset=\{x_1,...,x_N\}\subset \mathbb{R}^d$ is a point cloud.
We will denote the superlevel sets of $f$ by
\begin{equation}
    \label{eq:defsuplev}
    \levset(a)=f^{-1}[e^{-a},\infty) =
\left\{x\in \mathbb{R}^d:f(x)\geq e^{-a}\right\}
\end{equation}
using $a$ in negative log coordinates so that the family of subspaces is increasing.
Then the following is central problem in topology and data \cite{chazal2013bootstrap,phillips2015geometric,bubenik2012statistical,fasy2014confidence}:
\begin{prob}
\label{prob:superlev}
Determine the topological type (for instance, the persistent homology groups) of $\levset(a)$.
\end{prob}

The main step is to construct a \emph{filtered simplicial complex}, which is a simplicial complex $X$ together with a function $w:X\rightarrow \mathbb{R}$ defined on the simplices, with the property that $X(a)=w^{-1}(-\infty,a]$
is a subcomplex for all $a\in \mathbb{R}$.
In this setting, density levels are put on the same footing as the persistence values of higher order simplices, in contrast with a standard pipeline in Topological Data Analysis
that involves separately removing low-density data points as a preprocessing step.
In one key application \cite{chazal2011clustering}, the degree zero persistent homology of one such construction can be used as a rigorous form of density-based hierarchical clustering.
More generally, the superlevel sets of a kernel density estimator are a
topological model for the
density landscape of a point 
cloud with noise.

Constructing a suitable choice of 
$X(a)$ has proven to be challenging in higher dimensions.
A crucial but thorny step is to 
determine a set
$S \subset \mathbb{R}^d$ corresponding the vertices of $X$, which must be well-spaced, but also stable with respect to outliers or low-density regions.
To give a sense of some of the difficulties that arise, consider the following simple algorithm:
\begin{alg} \label{alg:intro}
Approximate the family $\levset(a)$ by a filtered simplicial complex $X(a)$.
\begin{enumerate}
    \item Choose a minimum separation value $r>0$, and initialize
    $S=\emptyset$.
    \item \label{item:naivealgsample} Sample 
    $M$ points $\dataset'=\{x'_1,...,x'_{M}\}$ from the underlying distribution of 
    $f(x)$. In other words, $x'_i$ is the sum of a uniformly random element of $\dataset$
    and a sample from the Gaussian kernel.
    \item 
    For each $x'_i$ in decreasing order of $f(x'_i)$, 
    make the replacement $S\mapsto S\cup \{x'_i\}$ if 
    $\lVert x_i'-p\rVert\geq r$
 for all existing points $p\in S$.
    \item Construct a filtered complex $(X,w)$ on the vertex set $S$,
    using general methods such as weighted Vietoris-Rips with $w(p)=-\log(f(p))$.
\end{enumerate}
\end{alg}
Algorithm \ref{alg:intro} seems intuitive, but it has undesirable propeties
that make it unusable in higher dimension. One is that
while the homotopy type of the superlevel sets 
$\levset(a)$ is a function only of the pairwise distances $\lVert x_i-x_j\rVert$, 
the expected size of the vertex set $S$ grows
rapidly by simply adding coordinates of 
zeroes to the end of every $x_i$, 
as it is proportional to 
the covering number of $\levset(a)$ by balls of radius $r$. 
Another issue is that the expected density of 
a sample point $f(x'_i)$, which is used to sort the points and filter homology, becomes small and 
dominated by noise for $d\gg 0$.

An alternate approach to producing a 
finite vertex set is
simply to discretize space,
but this is clearly only possible in very low dimensions.
Another is to enumerate 
the critical points of $f(x)$ and compute the discrete Morse complex \cite{forman1998morse,robins2011theory},
but this can in general lead to a combinatorial explosion in the number of points, including the 
case of sums of Gaussian kernels \cite{edelsbrunner2013risk}.
An intuitive fix is to use something similar to 
Algorithm \ref{alg:intro}, 
but using $\dataset$ itself in place of
$\dataset'$ in item \ref{item:naivealgsample}, so that the final answer is technically independent of the dimension of the embedding. However, 
this is not a stable solution 
as the size of $\dataset$ grows. Indeed, as
long as the underlying density from which $\dataset$ was sampled is nonzero on all of $\mathbb{R}^d$, a 
large enough sample will eventually fill out space,
so that we end up with the same poor scaling
with the covering number as in the original algorithm.

Once a density-weighted vertex set $S$ has been chosen, there are a number of persistence constructions 
that are robust with respect to noise, 
as well as statistically rigorous 
results about their output \cite{chazal2011geometric,phillips2015geometric,bobrowski2017topological,adcock2013ring,bubenik2012statistical,blumberg2012robust,carlsson2007multidimensional,carlsson2009computing,lesnick2011theory,blumberg2020stability}.

\subsection{Proposed method}

Instead of improving on Algorithm \ref{alg:intro}, our approach is to 
replace $\levset(a)$ itself with the superlevel sets 
of a sort of transformed function
$\tilde{f}:\mathbb{R}^d\rightarrow \mathbb{R}_+$,
given by
\begin{equation}
\label{eq:introftilde}
\tilde{f}(y)=\inf_{x\in \mathbb{R}^d}
f(x) \exp(\lVert x-y\rVert^2/2h^2).
\end{equation}
Let $\shape \subset \mathbb{R}^d$ be the region on which $\tilde{f}$ is nonzero, and let $\alpha:\shape\rightarrow \mathbb{R}$ be defined by $\alpha(x)=-\log(\tilde{f}(x))$. We then define
the \emph{shape} associate to $f(x)$
by
$\shape(a)=\tilde{f}^{-1} [e^{-a},\infty)=
\alpha^{-1} (-\infty,a]$.
By taking $x=y$, it is obvious that $\tilde{f}(y)\leq f(y)$, so that 
$\shape(a)\subset \levset(a)$, but not that either quantity retains interesting information. However, we will prove the following theorem.
\begin{thma}
We have
\begin{enumerate}
\item \label{item:introthmhull}
The full shape $\shape$ is the interior of the convex hull 
$\chull(\dataset)$, and $\alpha$ is a continuous function on it.
\item \label{item:introthmcover} We
have that
$\shape(a)\subset \levset(a)$, and
\[\levset(a)=
\bigcup_{y\in \shape(a)} B_r(y),\ 
r=\sqrt{2h^2(a-\alpha(y))},\] 
where $B_r(y)$ is the ball of radius $r$ centered at $y$.
    \item \label{item:introthmhomotopy}
    The inclusion map $i:\shape(a)\hookrightarrow \levset(a)$ induces a homotopy equivalence, whose inverse 
    is induced by the restriction of the map $p:\mathbb{R}^d\rightarrow \chull(\dataset)$ defined by
\begin{equation}
\label{eq:introconvmap}
    p(x)= \frac{\sum_{i}a_i \exp(-\lVert x-x_i\rVert^2/2h^2) x_i}{\sum_{i} 
    a_i
    \exp(-\lVert x-x_i\rVert^2/2h^2)}.
    \end{equation}    
\end{enumerate}
\end{thma}

Item \ref{item:introthmhomotopy} shows that using $\shape(a)$ in place of
$\levset(a)$ does not change the desired answer in Problem \ref{prob:superlev}.
On the other hand, the other parts show that for the purpose of computing a filtered complex $X(a)$, it is far better to select a vertex set $S$ from $\shape(a)$ rather than from 
$\levset(a)$. For one thing, item \ref{item:introthmcover} shows that selecting $S\subset \shape(a)$ corresponds to taking a subset of a cover, which in light of the nerve theorem
is far more natural than selecting vertices from the union.
On a practical level, the fact that $\shape(a)$ is a subset of the convex hull $\chull(\dataset)$ from item \ref{item:introthmhull} shows that 
the covering number of $\shape(a)$ does not have the poor
dependence on the embedding dimension $d$ that we had
with $\levset(a)$. It is not hard to see that the inverse homotopy equivalence $p(x)$ respects coordinate changes, and factors through the projection map onto the affine span of $\dataset$.

In Section \ref{sec:finplex},
we explain how to use $p(x)$ to replace the sampling step
in Algorithm \ref{alg:intro}, 
and it is not hard to see that this sampling procedure is
stable under coordinate changes as well. 
In order to generate filtered simplicial complexes, 
we will actually use the \emph{power diagram} in place of the covering by closed balls (which would be the weighted \v{C}ech cover), 
whose nerve is known as the alpha complex \cite{edelsbrunner1983shape}.
Alpha complexes have a number of 
theoretical benefits beyond persistent homology, namely that they are minimal in size, naturally embedded in space, and can be used to generate to beautiful geometric models \cite{edelsbrunner1983molecule}.
An illustration of the whole pipeline in spatial density estimation, including an image of $\shape(a)$, is shown in Figure
\ref{fig:species}.

In Section \ref{sec:experiments}, we present several experiments with this construction in higher dimension, resulting in compelling persistence homology calculations
as well as refined geometric models, taking advantage of the precise and visualizable natural of alpha complexes and alpha shapes. 
In order to make alpha complexes computable
in higher dimension, we used a recent scalable algorithm based on the 
duality principle \cite{carlsson2023alpha}, which has no explicit dependence on the embedding dimension. 
 All code for this project was written in MAPLE, and is available at the first author's website
 \href{http://www.math.ucdavis.edu/~ecarlsson/}{http://www.math.ucdavis.edu/~ecarlsson/}, along with supporting worksheets to generate each figure.

\subsection{Ackowledgments}
Both authors were supported by the Office of Naval Research, project (ONR) N00014-20-S-B001. 
    Thanks to Greg Kuperberg, who suggested using the Legendre transform to simplify an earlier version of Lemma 
\ref{lem:legendre}.

\section{Notation and preliminaries}

We summarize some preliminary definitions 
and notation for kernel density estimation and computational topology, including the power diagram and alpha complexes.

\subsection{Kernel density estimators}

\label{sec:kde}

Let $\dataset =\{x_1,...,x_N\}
\subset \mathbb{R}^d$ be a
point cloud. A Gaussian kernel density estimator is a sum of the form
\begin{equation}
\label{eq:kde}
f(x)=\sum_{i=1}^N a_i K_h(x-y),
\quad K_h(v)=\exp(-\lVert v\rVert^2/2),
\end{equation}
for $a_i>0$. 
We will be interested in the superlevel sets 
\begin{equation}
    \label{eq:suplev}
    \levset(a)=f^{-1}[e^{-a},\infty)=\left\{x:f(x)\geq e^{-a}\right\}
\end{equation}

For simplicity, we consider only finite sums, 
but our results apply to the convolution of more general distributions by Gaussian kernels. We will assume the the norm is always the standard $L^2$-norm, 
as any other quadratic form can 
be transformed in that way by a linear change of coordinates. Moving (anisotropic) metrics are an interesting extension, which we hope to study in future papers.
For now, we remark that general Riemannian metrics can often be approximated by Euclidean ones in higher dimensions, for instance using spectral embeddings, which will be used
in Section \ref{sec:experiments}. 
A more general setup might involve replacing
\eqref{eq:kde} by the convolution of a distribution on a Riemannian manifold by a heat kernel with respect to the metric.

\subsection{Computational topology}

A simplicial complex on a vertex set $S$
is a collection of $X=\{\sigma\}$
of nonempty subsets $\sigma \subset S$, which is closed
under taking nonempty subsets.
If $\{U_x:x\in S\}$ is a collection of subregions $U_x\subset \mathbb{R}^d$ indexed by $S$, then
the nerve is the simplicial complex on the vertex set $S$ defined by
\begin{equation}
    \label{eq:nerve}
    X=\nerve(\{U_x\}),
    \quad \{x_0,...,x_k\} \in X \Leftrightarrow
    U_{x_0} \cap \cdots \cap U_{x_k} \neq \emptyset.    
\end{equation}
One form of the \emph{nerve theorem} states
that if $U_x$ are convex, then the geometric realization is homotopy equivalent to the
union $|X| \sim \bigcup U_x$. For a reference on nerve theorems in general, see \cite{bauer2023unified}.

In persistent homology, a \emph{filtered} simplicial complex with vertices in $S$ is a nested family 
of complexes $X(a)$ for $a\in \mathbb{R}$, with the property that $X(a) \subset X(b)$ for $a\leq b$.
Equivalently, it is the data of a pair
$(X,w)$ where $X$ is a simplicial complex, and $w:X\rightarrow \mathbb{R}$ is a function with the property that the sublevel set 
$X(a)=w^{-1}(-\infty,a]$ is a complex for all $a$.
We then have the persistent homology groups
$H_*(X)$, which can be summarized by the barcode diagrams
\cite{Carlsson2009TopologyAD,carlsson2009computing,otter2015roadmap,edelsbrunner2010computational}. All barcode diagrams
generated for this paper were calculated
using \texttt{javaplex} \cite{adams2014javaplex}.

\subsection{Alpha complexes}

\label{sec:alpha}

For a reference on this section, we refer to 
\cite{aurenhammer1984optimal,
edelsbrunner1995union,
edelsbrunner2010computational}.
Let $S=\{p_1,...,p_n\}\subset 
\R^d$ be a collection of points, and let
$\pi :S\rightarrow \R$ be a function 
called the power map. We have the associated \emph{weight map} $w=w_{S,\pi}:\R^d\rightarrow \R$, defined by
\begin{equation}
\label{eq:powdist}
w(x)=\min_{p\in S} w_{p}(x),\quad w_{p}(x)=
\lVert x-p\rVert^2-\pi(p).
\end{equation}
Then the sublevel sets of $w$ can be expressed as a union of closed discs with varying radii
\begin{equation}
    \label{eq:weightunion}
    w^{-1} (-\infty,a]=\bigcup_{p\in S} D_{p}(a),
    \quad D_{p}(a)=\left\{x\in \mathbb{R}^d:w_p(x)\leq a\right\}.
\end{equation}
In other words, $D_p(a)=B_{p}(r)$ is a closed ball of radius $r=\sqrt{a-\pi(p)}$, or is empty if 
$a-\pi(p)<0$. Taking the nerve of the covering for each $a$ determines a filtered complex, known as the weighted \v{C}ech complex:
\begin{equation}
    X(a)=\left\{\{x_0,...,x_k\} \subset S: D_{x_0}(a) \cap \cdots \cap D_{x_k}(a) \neq \emptyset \right\}
\end{equation}

\begin{figure}
\centering
\begin{subfigure}[b]{.32\textwidth}
\includegraphics[scale=.2]{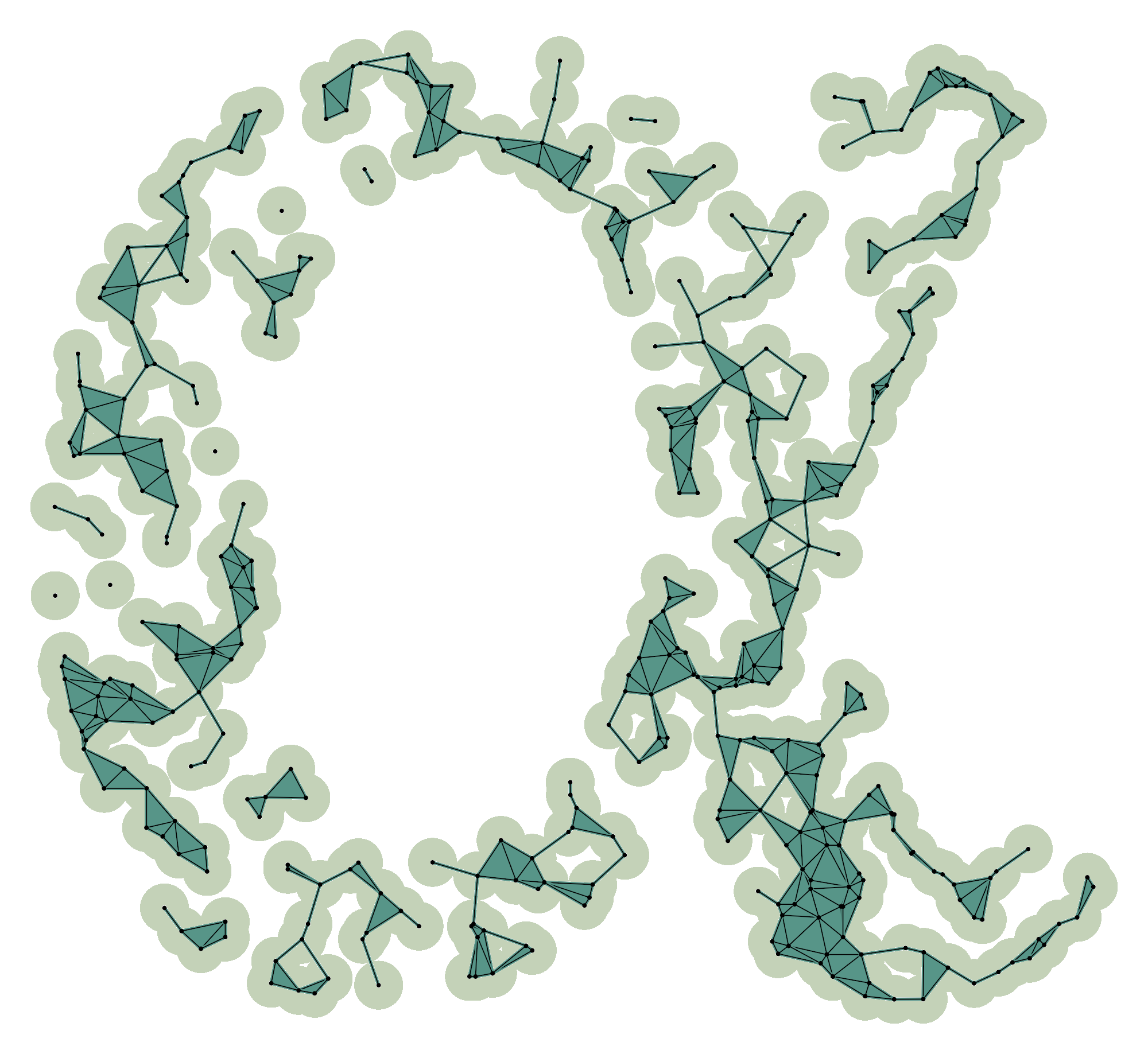}
\end{subfigure}
\begin{subfigure}[b]{.32\textwidth}
\includegraphics[scale=.2]{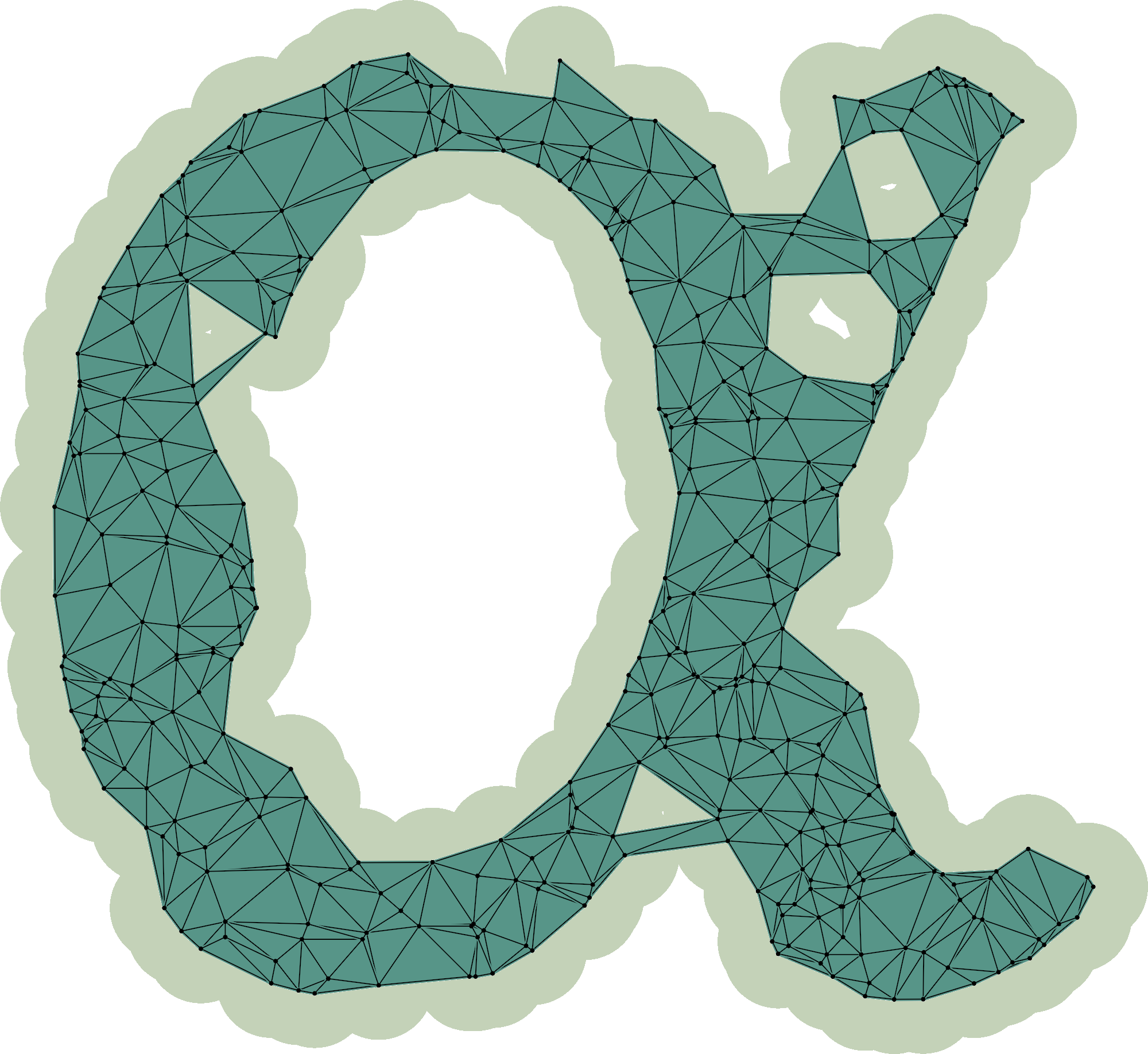}
\end{subfigure}
\begin{subfigure}[b]{.32\textwidth}
\includegraphics[scale=.2]{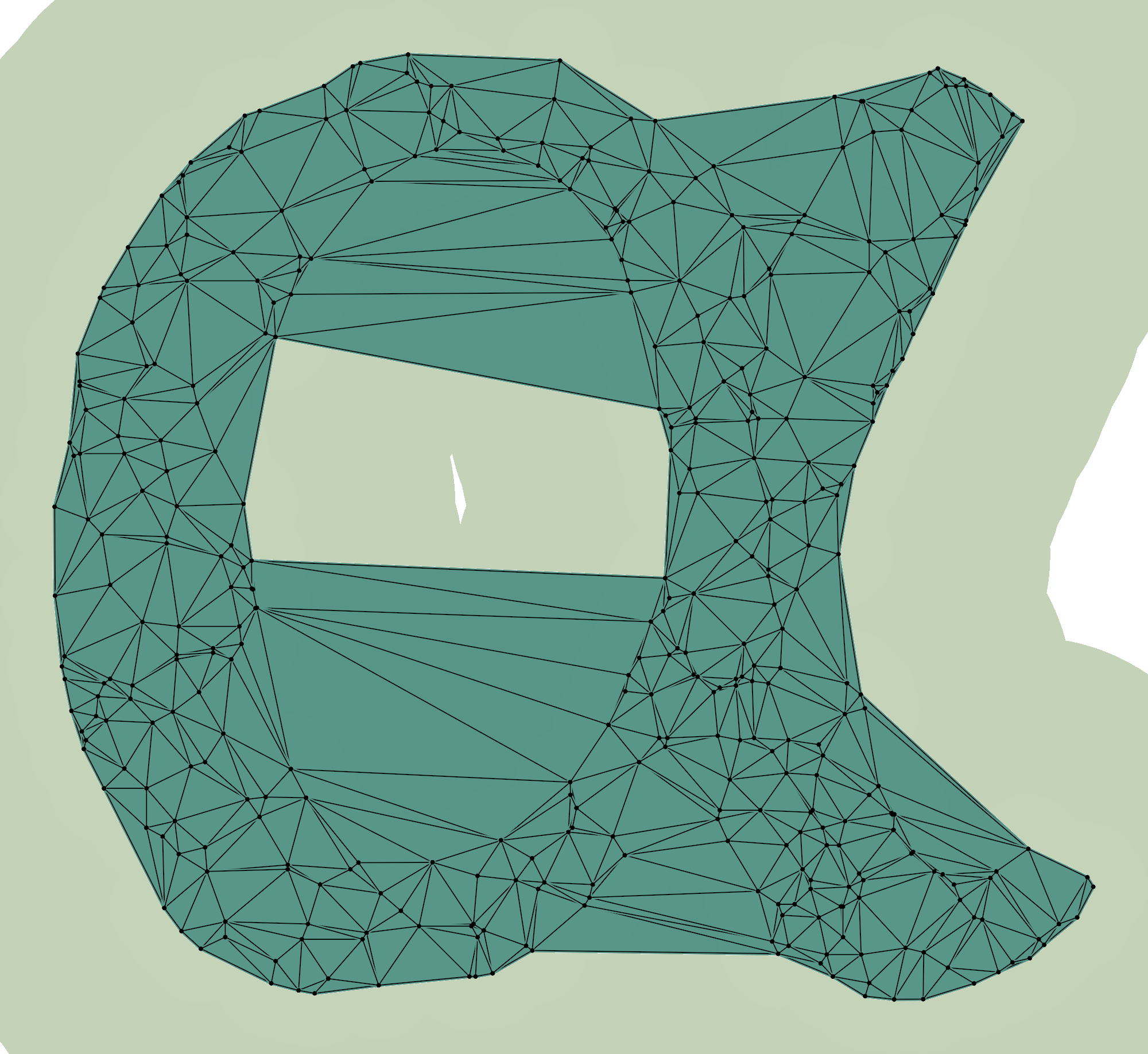}
\end{subfigure}
\caption{A sequence of power diagrams with their corresponding alpha shapes.}
\label{fig:alphalet}
\end{figure}

The power diagram $\{V_p(a)\}$ associated to 
$(S,\pi)$ comes from intersecting each element of the \v{C}ech cover with the weighted Voronoi cell
\begin{equation}
V_p(a)=D_p(a)\cap V_p,\quad 
    V_p=\left\{x \in \R^d :
w_{p}(x) \leq w_{p'}(x) \mbox{ for all $p'\in S$}\right\}
\end{equation}
The filtered complex $X(a)=\nerve(\{V_p(a)\})$ is known as the alpha complex. Since its union is the same as the union of discs, it is homotopy equivalent to the weighted 
\v{C}ech complex. Then geometric realization
\begin{equation}
\label{eq:alphashape}    
|X(a)|=\bigcup_{\sigma\in X(a)}
\chull(\sigma),
\end{equation}
where $\chull(\sigma)$ is the convex hull, is known as the alpha shape
\cite{edelsbrunner1983shape}.
In \cite{edelsbrunner1995union},
Edelsbrunner gave an explicit deformation retraction of the union of balls onto the shape. 
Figure \ref{fig:alphalet} illustrates the union of the cover and the corresponding shapes in the unweighted case.

\section{Alpha shapes of Gaussian KDE's}

We define our main constructions, including the proof of our main theorem, and our proposed algorithm for generating finite alpha complexes.
Let $\dataset=\{x_1,...,x_N\}\subset \mathbb{R}^d$, and consider a sum $f(x)$ of Gaussian kernels 
as in \eqref{eq:kde}. 
Let
\begin{equation}
\label{eq:ftilde}    
\tilde{f}(y)=\inf_{x\in \mathbb{R}^d}
f(x)K_h(x-y)^{-1}=
\inf_{x\in \mathbb{R}^d}
f(x)\exp(\lVert x-x_i\rVert^2/2h^2).
\end{equation}
We clearly have that 
$\tilde{f}(x)\leq f(x)$ by taking $y=x$
in the infimum. We will see in Lemma 
\ref{lem:legendre} that $f$ can be recovered from $\tilde{f}$, making it a sort of convex transform.
We will be interested in the corresponding superlevel sets of 
$\tilde{f}$.
\begin{defn}
\label{def:shape} 
Let $\shape=\{x\in \mathbb{R}^d:\tilde{f}(x)\neq 0\}$, and define $\alpha:\shape\rightarrow \mathbb{R}$ by $\alpha(y)=-\log(\tilde{f}(y))$.
The \emph{shape} associated to $f$ is the family of closed subregions
\begin{equation}
\label{eq:defshape}    
\shape(a)=\tilde{f}^{-1} [e^{-a},\infty)=
\alpha^{-1} (-\infty,a].
\end{equation}
\end{defn}

We now state our main theorem.
\begin{thm}
\label{thm:shape}
We have
\begin{enumerate}
\item \label{item:thmhull}
The full shape $\shape$ is the interior of the convex hull 
$\chull(\dataset)$, and $\alpha$ is a continuous function on it.
\item \label{item:thmcover} We
have that
$\shape(a)\subset \levset(a)$, and
\[\levset(a)=
\bigcup_{y\in \shape(a)} B_r(y),\ 
r=\sqrt{2h^2(a-\alpha(y))},\] 
where $B_r(y)$ is the ball of radius $r$ centered at $y$.
    \item \label{item:thmhomotopy}
    The inclusion map $i:\shape(a)\hookrightarrow \levset(a)$ induces a homotopy equivalence, whose inverse 
    is induced by the restriction of the map $p:\mathbb{R}^d\rightarrow \chull(\dataset)$ defined by
\begin{equation}
\label{eq:convmap}
    p(x)= \frac{\sum_{i}a_i \exp(-\lVert x-x_i\rVert^2/2h^2) x_i}{\sum_{i} 
    a_i
    \exp(-\lVert x-x_i\rVert^2/2h^2)}.
    \end{equation}    
\end{enumerate}
\end{thm}

\subsection{Proof of Theorem \ref{thm:shape}}

We now prove Theorem \ref{thm:shape} in a sequence of lemmas. 

\begin{lemma}
\label{lem:legendre}    
We have
\begin{equation} \label{eq:proplegendre} f(x)=\sup_{y\in \mathbb{R}^d}
\tilde{f}(y) \exp(-\lVert x-y\rVert^2/2h^2).
\end{equation}
Moreover, for each $x$, the supremum in \eqref{eq:proplegendre} is obtained at the unique 
value $y=p(x)$, where $p(x)$ is the function from \eqref{eq:convmap}.

\end{lemma}

\begin{figure}
    \centering
    \includegraphics[scale=1.2]{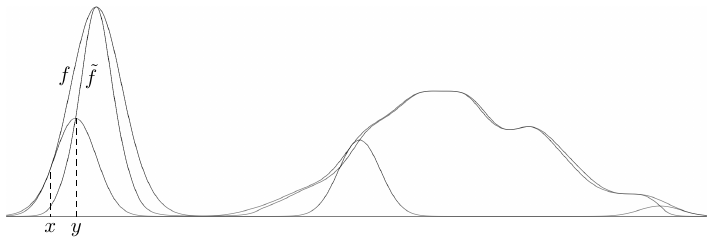}
    \caption{A univariate Gaussian KDE $f$ shown as the top curve, with $\tilde{f}$ shown as the bottom curve. The Gaussians $g(x)$ all have the same scale as was used to define $f(x)$, but different scalar multipliers and centers. Lemma
    \ref{lem:legendre} shows that $f$ and $
\tilde{f}$ determine each other, and that $y=p(x)$
    for any pair $(x,y)$.}
    \label{fig:univariate}
\end{figure}

An illustration is shown in Figure
\ref{fig:univariate}.

\begin{proof}

First, we claim that the function
\begin{equation}
    \label{eq:convfunc}
F(x)= \lVert x\rVert^2/2+h^2\log(f(x))
\end{equation}
is convex.
Then we may write
$\alpha(y)=\lVert y \rVert^2/2h^2-F^*(y)/h^2$,
 where 
\begin{equation}
\label{eq:legendre}
    F^*(y)=\sup_{x} \left(x\cdot y-F(x)\right)
\end{equation}
is the Legendre transform.
The first statement then follows from the Fenchel-Moreau theorem, which states
that $F^{**}=F$.

To prove the convexity, we start with a useful fact: suppose $\mathbb{A}\subset \mathbb{R}^d$ is an affine subspace, 
and let $f|_{\mathbb{A}}$ be the restriction of $f$ to $\mathbb{A}$.
Then
    \begin{equation}
        \label{eq:affkde}
        f|_{\mathbb{A}}(x)=\sum_{i} b_i \exp(-\lVert x-\overline{x}_i\rVert^2/2h^2),\quad
        b_i=a_i \exp(-\lVert x_i-\overline{x}_i \rVert^2/2h^2).
    \end{equation} 
where $\overline{x}_i$ is the orthogonal projection of $x_i$ onto $\mathbb{A}$.
Then it suffices to prove
the convexity of a one-dimensional function $F(x+tv)$,
which by \eqref{eq:affkde} reduces to 
the one-dimensional case. For this, it suffices to show that the second derivative is nonnegative. 
We check
\[f''=h^{-2}\left(\sum_{i} 
a_i(x-x_i)^2 e^{-(x-x_i)^2/2h^2}\right)
-f.\]
We can now write
\[F''=\frac{f''}{f}-\left(\frac{f'}{f}\right)^2+1=
h^{-2}\left(
E[(X_i-x)^2]-(E[X_i-x])^2\right)\]
where the expectations are with respect to the probability measure on $\dataset$ in which each $x_i$ 
is weighted in proportion to
$\exp(-\lVert x-x_i\rVert^2/2h^2)$.
Noting that this is the expression for the variance, 
we obtain the desired nonnegativity.

For the second statement, suppose $y$ is any value that obtains the supremum in \eqref{eq:proplegendre} for a given $x$,
and let $g(z)=c \exp(-\lVert z-y\rVert^2/2h^2)$ be the corresponding Gaussian centered at $y$ with
$c=\tilde{f}(y)$ as in Figure \ref{fig:univariate}.
Then we have that $g$ agrees with 
$f$ to first order at $x$:
\begin{equation}
\label{eq:lemfit}
g(x)=f(x),\quad 
\grad_{g} (x)=\grad_{f} (x).
\end{equation}
Dividing the
second second expression by the first
and solving for $y$, we obtain
the expression in \eqref{eq:convmap}, 
which also establishes 
the uniqueness.
\end{proof}

Taking the log in Lemma \ref{lem:legendre}, we obtain the following useful formula:
\begin{equation}
\label{eq:alphaofp}    
\alpha(p(x))=-\log(f(x))-
\lVert x-p(x)\rVert^2/2h^2.
\end{equation}
The next lemma shows that every point of $\shape$
is in the image of $p$, so that \eqref{eq:alphaofp} gives a complete description of $\shape(a)$.

\begin{lemma}
    \label{lem:hull}
    The image of $p$ is equal to the interior of the convex hull $\chull(\dataset)$, which is equal to 
the full shape $\shape$. The fiber of $p$ over 
$p(x)$ is the perpendicular affine subspace to 
$\affhull(\dataset)$
passing through $x$.
\end{lemma}

\begin{proof}

First, if $p|_{\mathbb{A}}$ is the restriction to an affine subspace, we have
\begin{equation}
\label{eq:projp}
    \Pi(p|_{\mathbb{A}}(x))=
    \frac{\sum_i b_i \exp(-\lVert x-\overline{x}_i\rVert^2/2h^2)\overline{x}_i}
    {\sum_i b_i\exp(-\lVert x-\overline{x}_i\rVert^2/2h^2)}
\end{equation}
where $\Pi$ is the orthogonal projection, and $b_i$ and $\overline{x}_i$ are as in \eqref{eq:affkde}.
Then for the second statement,
it suffices to show that if 
$\affhull(\dataset)=\mathbb{R}^d$ then 
$p$ is injective.
In this case, suppose that 
$p(x)=p(x')$ for $x \neq x' \in \R^d$.
Then applying \eqref{eq:projp} again to 
$\{x+tv\}\subset \mathbb{R}^d$ with $v=x'-x$ reduces the problem to the one dimensional case, noting that the condition that $\dataset$ affinely spans $\mathbb{R}^d$ means that the projected data set cannot all map to the same point. The 
one-dimensional case can 
then easily be checked by
showing that $p'(t)$ is positive,
so that $p$ is increasing, which is an easy calculation using the quotient rule.

For the first statement,
suppose that $x$ is not in $\chull(\dataset)$, or is on the boundary. Let $v$ be the outward-pointing unit normal vector to any face of 
$\chull(\dataset)$ for which $x$ is not strictly on the interior side. We can then check that $\tilde{f}(x)=0$
by setting $y=x+tv$ in \eqref{eq:introftilde} and taking the limit as $t\rightarrow \infty$,
showing that $\shape$ is contained in the interior.
Now by \eqref{eq:alphaofp}, we have that $p(\mathbb{R}^d) \subset \shape$, so it remains to show that
every point $y$ in the interior is equal to some $p(x)$. 
Using \eqref{eq:projp} again, we may assume that $\affhull(\dataset)=\mathbb{R}^d$. 

We first show that if $y$ is any point
on the boundary of $\chull(\dataset)$, then $y$ is a limit of points in the image of $p$.
Let $\mathbb{A}$ be the affine subspace 
of codimension $1$ containing any face that contains $y$, and let $v$ be the outward-pointing unit normal vector. Then we have that
\[\lim_{t\rightarrow \infty} p(y+tv)=
\overline{p}(y),\]
where $\overline{p}:\mathbb{A}\rightarrow \mathbb{A}$ is the result of setting
$a_i=0$ for $x_i \notin \mathbb{A}$
in \eqref{eq:convmap}. By induction on the dimension, we find that the boundary of $\chull(\dataset)$ is in the closure of the image of $p$.
Since the image $p(\mathbb{R}^d)$ is contractible, it must contain the interior of $\chull(\dataset)$.

\end{proof}

\begin{lemma}
\label{lem:retract}    
We have that $p(\levset(a))\subset \shape(a)$ is a deformation retract.
\end{lemma}

\begin{remark}
In \cite{edelsbrunner1995union}, Edelsbrunner gave an explicit deformation retraction of the union of balls onto the shape. However,
in our deformation retraction from Lemma \ref{lem:retract}, $\shape(a)$ is the outer set.
\end{remark}

\begin{proof}

If the affine span $\mathbb{A}=\affhull(\dataset)$ is a strict subset of $\mathbb{R}^d$, 
then we can easily see that $\levset(a)$ retracts
onto $\levset(a)\cap \mathbb{A}$.
By Lemma \ref{lem:hull}, we have that
$\shape(a)\subset \mathbb{A}^d$, and
using \eqref{eq:projp}, we may assume that 
$\affhull(\dataset)=\mathbb{R}^d$. In this case the restriction of $p$ 
to $\levset(a)$ is a homeomorphism onto its image $p(\levset(a))$, because the inverse of a continuous bijection on a compact set is continuous.
It suffices to produce a 
deformation retraction of
$\levset(a)\subset p^{-1}(\shape(a))$. 

Define $\varphi:p^{-1}(\shape(a)) \rightarrow \levset(a)$ by 
setting $\varphi(x)$ to be the point on $B_{r}(p(x))$ which is closest to $x$, with $r=\sqrt{2h^2(a-\alpha(p(x)))}$.
The statement that the image is contained in 
$\levset(a)$ follows from
\eqref{eq:alphaofp}. It also follows 
from \eqref{eq:alphaofp} that $\varphi$ acts as the identity on $\levset(a)$, since each $x\in \levset(a)$ is in $B_r(p(x))$.
Note that $\varphi(x)$ is always on the line segment connecting $x$ with $p(x)$.
We claim that the function
defined by $F(x,t)=(1-t)x+t\varphi(x)$ is the desired deformation
retraction. All that needs to be checked 
is that when $x\in p^{-1}(\shape(a))$,
we have $F(x,t) \in p^{-1}(\shape(a))$ for all $t$. In other words, we must check that for any $t\in [0,1]$,
we have $\alpha((1-t)x+tp(x))\leq a$.

We first claim that the expression
\begin{equation}
\label{eq:defleq}
A(x')=-\log(f(x'))-\lVert x'-p(x)\rVert^2/2h^2.
\end{equation}
is maximized at $x'=x$, which can be checked by showing the gradient is zero in every direction, and that noting that $A(x')$ 
is the negative of a 
convex function (in fact for any point in place of $p(x)$) by the proof of Lemma \ref{lem:legendre}.
Since $A(x)=\alpha(p(x))$, we have that $A(x')$ is bounded above by $a$ for all $x'$. Then it suffices to check that
\begin{equation}
    \label{eq:defdiff}
   \alpha(p(x'))\leq A(x') \Longleftrightarrow
    \lVert x'-p(x')\rVert^2 \geq
    \lVert x'-p(x)\rVert^2
\end{equation}
for $x'=(1-t)x+tp(x)$ with $t\in [0,1]$, 
so that $\alpha(p(x'))\leq a$.
To do this, we notice that the projection of $p(x')$ onto the line $\mathbb{A}=\{x+tv\}$ for $v=x'-x$ can only be closer
to $x'$ than $p(x')$. Then
\eqref{eq:projp} reduces the right hand side of \eqref{eq:defdiff} to the one-dimensional-case. The one-dimensional case follows because $p$ is increasing by the proof of Lemma \ref{lem:hull}, and because $p(x')$ is always on the opposite side of $p(x)$ from $x'$ for $t\in [0,1]$.

\end{proof}

\begin{proof}[Proof of Theorem \ref{thm:shape}]

Item \ref{item:introthmhull}
follows from Lemma \ref{lem:hull}
and \eqref{eq:alphaofp}. Item
\ref{item:introthmcover} follows from taking the log in Lemma \ref{lem:legendre}.
For item \ref{item:introthmhomotopy}, we have that $p(\levset(a))$ is homeomorphic to $\levset(a) \cap \affhull(\dataset)$, 
which is homotopy equivalent to $\levset(a)$.
By Lemma \ref{lem:retract}, we then have that
$\levset(a)$ is homotopy equivalent to 
$\shape(a)$. This also shows the equivalence is induced by the inclusion $\shape(a)\subset \levset(a)$, and that the inverse in induced by $p$.

\end{proof}

\subsection{Algorithm for generating finite alpha complexes}

\label{sec:finplex}

We now have the following modification of Algorithm \ref{alg:intro}, using the shape
$\shape(a)$ in place of $\levset(a)$.

\begin{alg} \label{alg:densland}
Represent $\shape(a)$ by a weighted alpha complex.
\begin{enumerate}
\item Initialize $S=\emptyset$. Fix a number $0<s<1$, and let 
$\epsilon=-\log(s)$.
    \item \label{item:algsample} Sample $M$ points $\{x'_1,...,x'_{M}\}$ 
    from the underlying distribution of $f(x)$, 
    and let $\dataset'=\{y_i'\}$
    where $y'_i=p(x'_i)$.    
    \item Scanning over all 
    $y'_i$ in decreasing order of $\tilde{f}(y'_i)$, make the replacement $S\mapsto S\cup \{y'_i\}$ if $y'_i$ has squared distance at least
    $2h^2\epsilon$ from all 
    existing points of $S$.
    \item Return the weighted alpha complex 
$X(a)=\nerve(\{V_p(2h^2a)\})$
  associated to the set $S$,
    with $\pi(p)=-2h^2\alpha(p)$
\end{enumerate}
\end{alg}

\begin{figure}[ht]
\centering
\begin{subfigure}[b]{.49\textwidth}
\centering
\includegraphics[scale=.9]{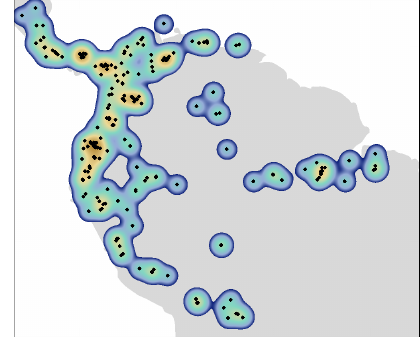}
\caption{$-\log(f(x))$ restricted to $\levset(a_0)$.}
\end{subfigure}
\begin{subfigure}[b]{.49\textwidth}
\centering
\includegraphics[scale=.9]{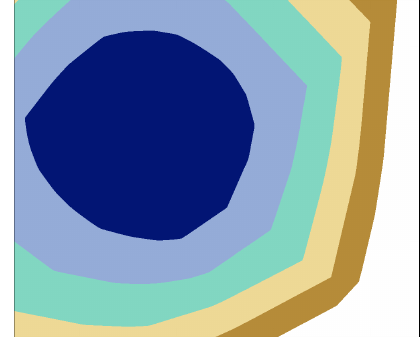}
\caption{$F(x)=\lVert x\rVert^2/2+h^2\log(f(x))$.}
\end{subfigure}
\centering
\begin{subfigure}[b]{.49\textwidth}
\centering
\includegraphics[scale=.9]{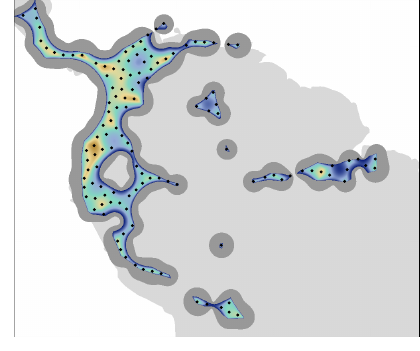}
\caption{$\alpha(y)=\lVert y\rVert^2/2h^2-F^*(y)/h^2$, restricted to 
$\shape(a_0)$.}
\end{subfigure}
\begin{subfigure}[b]{.49\textwidth}
\centering
\includegraphics[scale=.9]{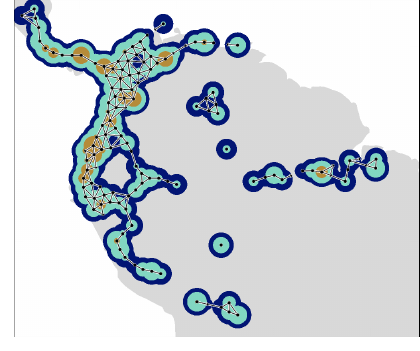}
\caption{Alpha complex
$\densalpha(f,\dataset',.3,d_0)$}
\end{subfigure}
\caption{Top left: a heat map of a kernel
density estimator 
$f(x)$ for the a well-known data set 
from spatial density estimation
\cite{phillips2006maximum}, cut off below a particular minimum density value $d_0=e^{-a_0}$.
Top right: The convex function $F(x)$, used in the proof of Lemma \ref{lem:legendre}.
Bottom left: a similar heat map of $\alpha$
restricted to $\shape(a_0)$
with $\levset(a_0)$ shown in gray, illustrating the homotopy equivalence in 
Theorem \ref{thm:shape}. 
To display $\alpha(y)$, we rasterized $F$ to give a $961\times 1061$ image, and used
the linear time Legendre transform
\cite{lucet1997faster} to evaluate $F^*$.
A set $S$ of landmark points resulting from
Algorithm \ref{alg:densland} with $s=.3$ and a density cutoff of $d_0$ is shown on top. Lower right: the resulting power diagram and its corresponding alpha complex.
}
\label{fig:species}
\end{figure}

An illustration is shown in the case of spatial density estimation in Figure \ref{fig:species}. 
We have chosen to represent the minimum spacing in terms of a unitless parameter $s$. For instance, a value of $s=.9$ would always correspond fine spacing, whereas $s=.3$ would have points fairly far apart. By properties of the Gaussian and the explicit form of \eqref{eq:convmap}, it is not hard to see that the sampling algorithm itself is not affected by the addition of new coordinates. In other words, the pushforward measure to $\shape(a)$ under $p$ is respected by coordinate changes.
Note also that 
$\alpha(y'_i)$ is readily computable
using \eqref{eq:alphaofp}.

\begin{defn}
The complex resulting from Algorithm \ref{alg:densland}
will be denoted $X=\densalpha(f,\dataset',s)$.
We will also denote $\densalpha(f,\dataset',s,d_0)=
\densalpha(f,\dataset'',s)$ where $\dataset''=\{y'_i\in \dataset':\tilde{f}(y_i')\geq d_0\}$.
\end{defn}

One artifact of having well spaced points in $S$ is that the spacing introduces additional connected components. In terms of persistence, this creates many 
intervals in 
the $\betti_0$ barcode diagram of length proportional to 
$\epsilon$. A reasonable way to 
compensate for this is to define
$\densalphao(f,\dataset',s)$ to be the result of replacing
\begin{equation}
    \label{eq:shiftvariant}
    w(\sigma)\mapsto \max(w(\sigma),\alpha(p_0)+\epsilon,...,\alpha(p_k)+\epsilon)
    \end{equation}
    for all simplices in $\densalpha(f,\dataset',s)$,
where $\epsilon=-\log(s)$ is the spacing parameter from 
Algorithm \ref{alg:densland}. This is the same as having each element of
$\{V_p(2h^2a)\}$ appear once its radius has reached the reached the minimum spacing threshold.

\section{Experiments}

We constructed the filtered complexes from Section \ref{sec:finplex} in several examples. 
To compute alpha complexes in higher dimension, we used a recent algorithm based on dual programming \cite{carlsson2023alpha}.
All alpha complexes in this paper took at most a few seconds to calculate, whereas the higher dimensional examples would have been impossible using other methods, because nearly all of them begin by computing the full (weighted) Delaunay triangulation.
In the sampling step we usually chose 
$M=|\dataset'|$ to be a few thousand,
which often took a couple of minutes,
coming from evaluating $p$ at that many points.
We will also write
$d_0=\quantile_f(p)$ to denote the quantile rank of $0\leq p\leq 1$, which is greatest density level at which the fraction of $x_i\in \dataset$ for which $f(x_i)\geq d_0$ is at least $p$.
A MAPLE worksheet showing how each example was done is available at the first author's website, under software: 
 \href{http://www.math.ucdavis.edu/~ecarlsson/}{http://www.math.ucdavis.edu/~ecarlsson/}.

\label{sec:experiments}

\subsection{Interesting energy landscapes}

\label{sec:config}

\begin{figure}[ht]
    \centering
    \begin{subfigure}[b]{.9\textwidth}
        \centering
        \includegraphics[scale=1]{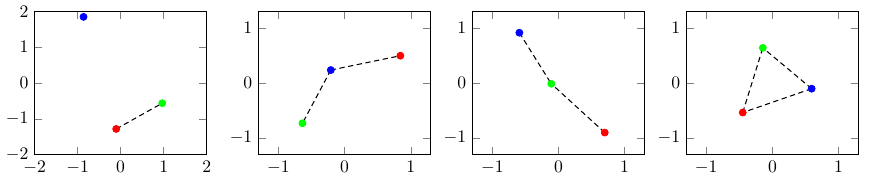}
    \end{subfigure}    
    \caption{Some
        mean-centered point configurations from the particle configuration model, using dashed lines to indicate a distance of roughly $1$.}
        \label{fig:configpoints}
    \end{figure}

We considered a simple energy function for three particles in the plane
\begin{equation}
    \label{eq:configloss}
H(p_1,p_2,p_3)=
\sum_{i<j} V(\lVert p_i-p_j\rVert)
+c\sum_{i} \lVert p_i\rVert^2
\end{equation}
where $V(r)=4(r^{-12}-r^{-6})$ is the
Lennard-Jones potential, which rewards pairs of particles that are approximately distance 1 apart, strongly penalizes particles that are much closer than that, 
and is neutral for far away points. 
The second factor in \eqref{eq:configloss} is there to keep points from wandering off to infinity.
Some typical configurations are shown in Figure \ref{fig:configpoints}.

\begin{figure}
    \centering
    \begin{subfigure}[b]{.25\textwidth}
        \centering
        \includegraphics[scale=.14]{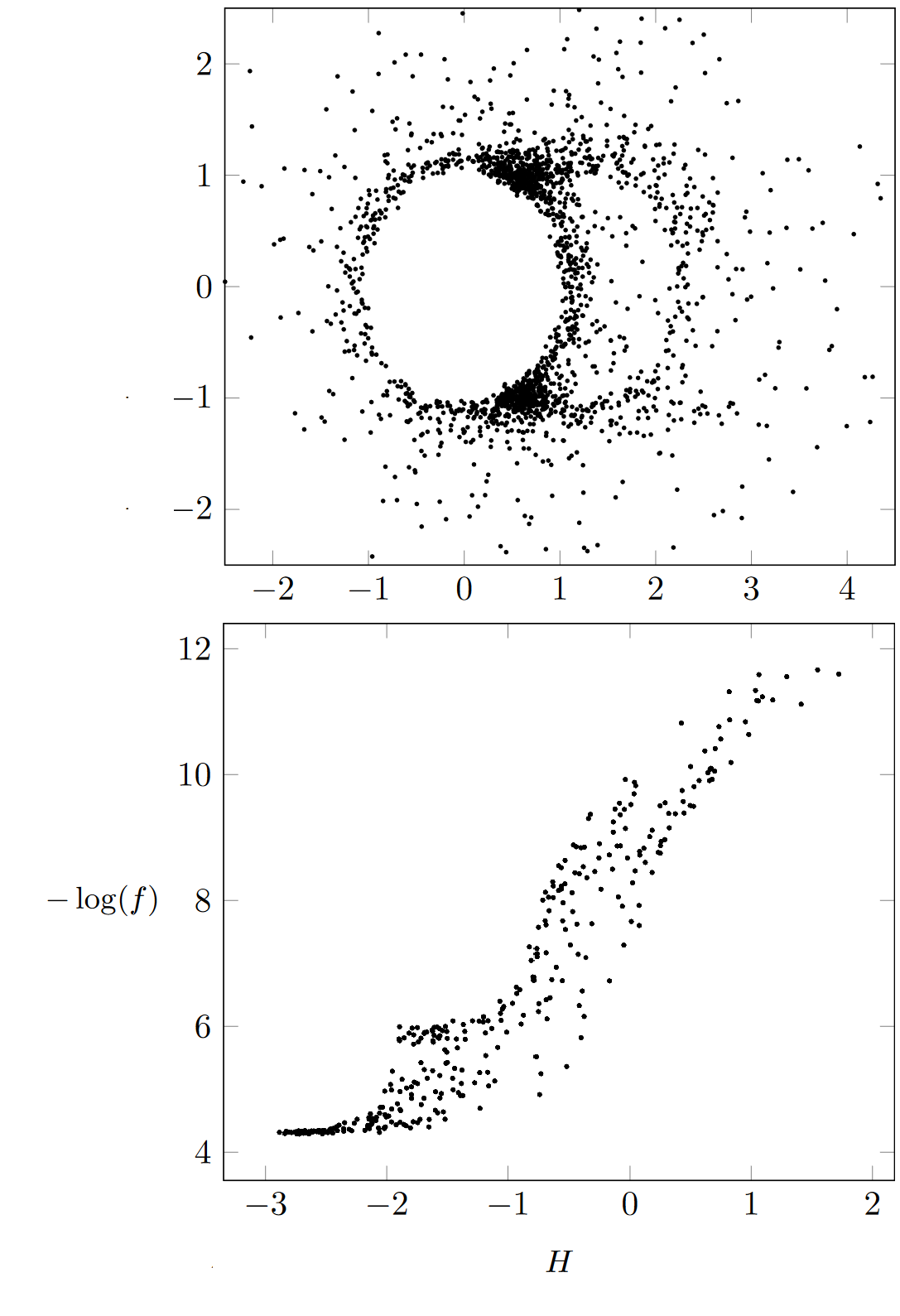}
    \end{subfigure}   
    \begin{subfigure}[b]{.6\textwidth}
    \centering
    \includegraphics[scale=1.45]{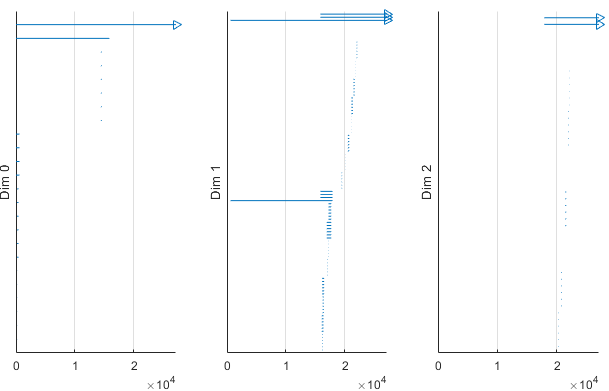}   
    \end{subfigure}
        \caption{Top left: a 2-dimensional point cloud of particle configurations
        obtained by normalized the first and second points to be at $(0,0)$ and $(1,0)$, and plotting the third.
        Bottom left: a plot of 
        $-\log(f(x_i))$ versus energy value, using $h=.3$.
        Right:
        \texttt{javaplex} output for 
        $\densalphao(f,\dataset',s,d_0)$
        using $M=60000$ samples, $s=.7$, and $d_0=\quantile_f(.6)$. }
    \label{fig:configbetti}     
\end{figure}

Using the Metropolis algorithm, we sampled 120000 points from the distribution
\[\rho(p_1,p_2,p_3)=
\exp(-\beta 
H(p_1,p_2,p_3))\] 
with $\beta=3.0$,
and mean-centered the output, producing a point cloud $\mathcal{D}$ in a 4-dimensional subspace of $\mathbb{R}^6$. We then formed the corresponding density estimator $f(x)$ with $h=.3$. 
In Figure \ref{fig:configbetti}, we see a normalized point cloud of configurations, and a plot showing a reasonable relationship between kernel density and energy.


We computed the alpha complex
$\densalphao(f,\dataset',s,d_0)$
using $s=.7$ and $d_0=\quantile_{f}(.6)$ up to the $3$-simplices, resulting in the sizes 
$(|X_0|,|X_1|,|X_2|,|X_3|)=
(2298,24582,64896,64674)$.
The persistent homology groups,
shown on the right side of Figure
\ref{fig:configbetti}, are consistent with a disjoint union of two circles turning into a configuration space of 3 ordered points in the plane.
This is because the densest points are the ones forming an equilateral triangle, which come in two types corresponding to he rotationally inequivalent permutations of the labels.
The points with two connections form the shape of a configuration space, which has betti numbers of 
$(\betti_0,\betti_1,\betti_2)=(1,3,2)$ corresponding to the arrows in the diagram \cite{cohen1995configuration}.
The fact that the barcodes have hardly any noise has to do with the 
precise nature of the alpha complex.

\subsection{Multiple optimizers in nonlinear regressions}

We next use alpha complexes to demonstrate the presence of 
local basins in the loss function of a nonlinear regression, using the Metropolis algorithm and maximum likelihood estimation to generate a point cloud. 
For our application we have chosen Example 2.6 from \cite{seber2005nonlinear},
which deals with the catalytic isometrization of $n$-pentane to $i$-pentane in the presence 
of hydrone, based on an original study by
Carr \cite{carr1960kinetics}.
The training data consists of
24 experimental runs, with four columns labeled $x_1,x_2,x_3,r$,
measuring the partial pressures of
hydrogen, $n$-pentane, and $i$-pentane, and the corresponding reaction rate $r$.
The modeling problem is to predict the reaction rate from the other variables using the model
\begin{equation}
    \label{eq:reaction}
r \sim
\frac{\theta_1\theta_3(x_2-x_3/1.632)}
{1+\theta_2 x_1+\theta_3 x_2+\theta_4 x_3}.
\end{equation}

\begin{figure}
\centering
        \begin{subfigure}[b]{.3\textwidth}
        \centering
        \includegraphics[scale=.5]{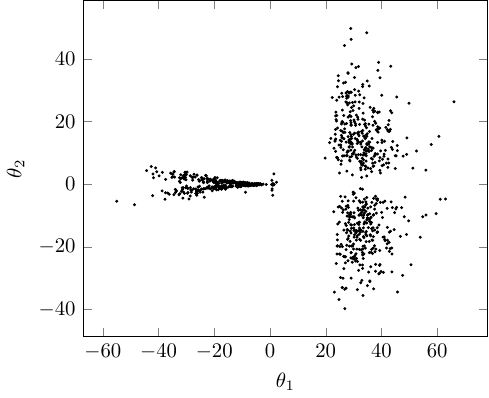}
    \end{subfigure}
    \begin{subfigure}[b]{.3\textwidth}
    \centering
        \includegraphics[scale=.515]{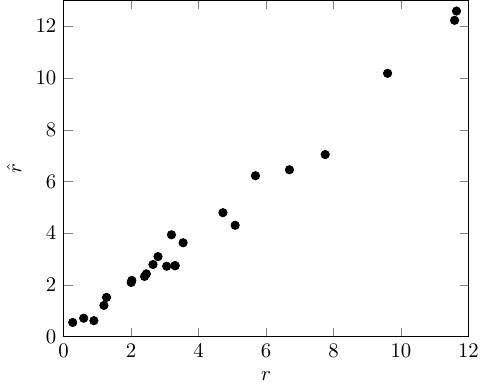}
    \end{subfigure}
        \begin{subfigure}[b]{.3\textwidth}
    \centering
        \includegraphics[scale=.515]{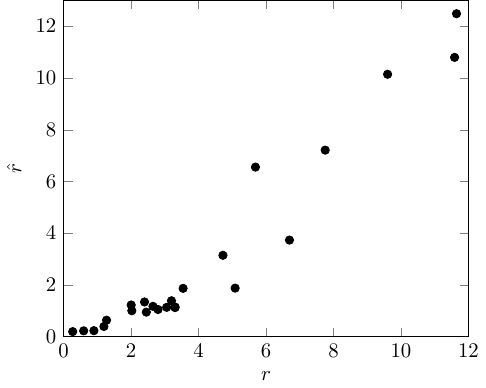}
    \end{subfigure}
    \caption{Left: a scatter plot of the first two parameters $(\theta_1,\theta_2)$ resulting from modeling the Carr data. In the middle/right frames we have the plot of predicted value versus true value for a typical value from the positive/negative values of $\theta_1$ groups respectively. The $\theta_1>0$ group appears to produce better fits.}
    \label{fig:carr}
\end{figure}

We simulated 10000 samples of the 
$\theta$-parameters using the Metropolis
algorithm, and a maximum likelihood with squared residual losses, 
\begin{equation}
    \label{eq:carrloss}
    p(\theta)=\exp(-\beta L(\theta)),\quad 
    L(\theta)=
    \sum_{i=1}^{24} \left(r_i-\hat{r}_i(\theta)\right)^2
\end{equation}
with a temperature of $\beta=3.0$,
resulting in a point cloud $\dataset \subset \mathbb{R}^4$.
One finds that there was a near-symmetry in simulaneously sending 
$\theta_i\mapsto -\theta_i$ 
for $i=\{2,3,4\}$, due to the fact that the parameters tended to dominate the leading
$1$ in the denominator.
Additionally, there appear to be local basins 
of solutions corresponding to the sign of 
$\theta_1$, with the ``correct'' positive values corresponding to better predictions. The results are shown in Figure \ref{fig:carr}.

An obvious next step would be to apply kernel density estimation to the point cloud 
$\dataset \subset \mathbb{R}^4$,
with some choice of scale $h$, but this would not be meaningful, as it is 
depends on the parametrization of the model. Instead, we define $f(x)$ using a new point cloud $\dataset\subset \mathbb{R}^{24}$ by mapping each $\theta$ to the corresponding vector of predictions
\begin{equation}
\label{eq:predmap}
\theta\mapsto (\hat{r}_1(\theta),...,\hat{r}_{24}(\theta)).
\end{equation}
We then have a natural choice of the scale parameter, $h=(\beta/2)^{-1/2}\sim .816$.

\begin{figure}
\centering
\begin{subfigure}[b]{1\textwidth}
\centering
     \includegraphics[scale=.55]{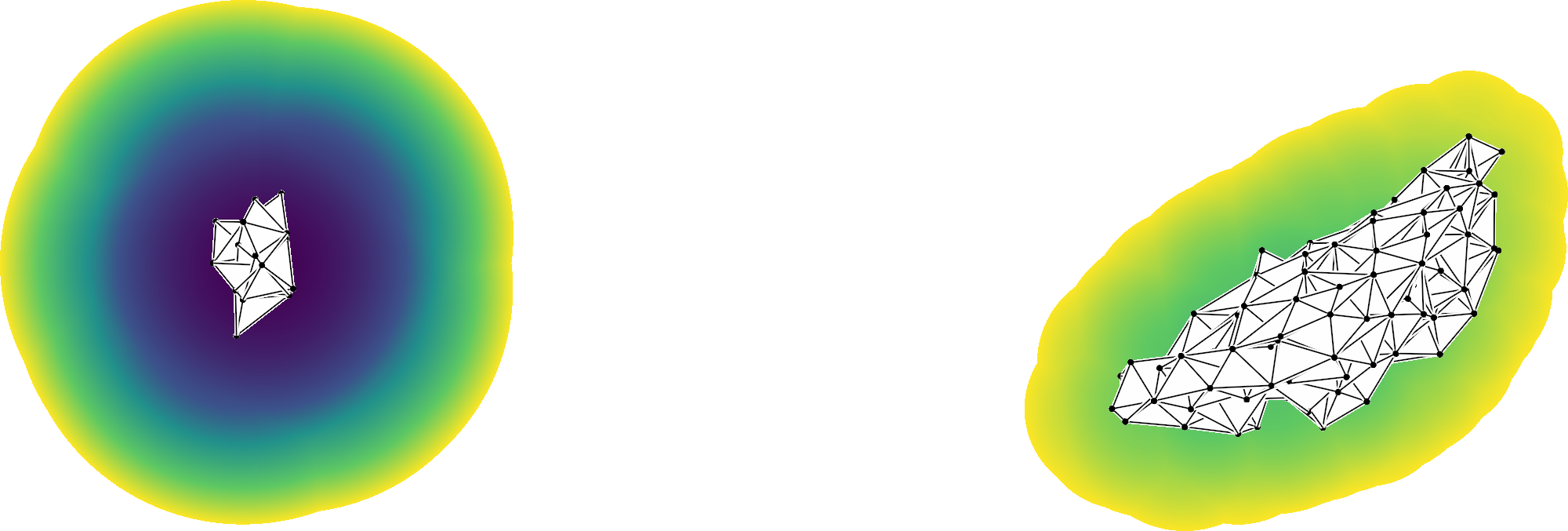}     
     \caption{Alpha complex projected to 2 dimensions}
\end{subfigure}
  \begin{subfigure}[b]{1\textwidth}
    \centering
    \includegraphics[scale=2]{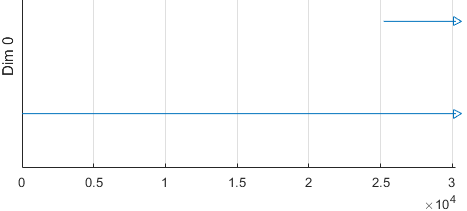}    
     \caption{Zero-dimensional persistence barcodes}
\end{subfigure}
    \caption{Alpha complex $\densalphao(f,\dataset,.6)$ for samples from maximum likelihood estimation 
    for the Carr model. The denser region on the left is the $\theta_1>0$ group, and the right cluster is $\theta_1<0$. The complex is calculated in $\mathbb{R}^{24}$, and then orthogonally projected for display purposes onto the first two principle components of
    $\dataset \subset \mathbb{R}^{24}$. The heat map represents the power diagram, which is also projected onto those coordinates. In the second frame, we have the $\betti_0$ barcode diagram, exhibiting the local basins.}
    \label{fig:carralpha}
\end{figure}

We built the alpha complex using a value of 
$s=.6$, and density cutoff $d_0=\quantile_f(.8)$, shown in Figure \ref{fig:carralpha}. 
The map 
to $\mathbb{R}^{24}$ effectively collapses the symmetry arising from the sign changes, leaving only
the two connected components. The component corresponding
to positive values of $\theta_1$ results in more accurate predictions, resulting in the denser component on the
left. The others lead to local basin which wider but less dense, corresponding to higher values of the loss function.

\subsection{A simple singular model}

We apply a similar method from the previous example to a singular statistical model,
which is Example 1.2 from \cite{watanabe2009algebraic}, 
with thanks to Dan Murfet for the suggestion. We see that the alpha complex exhibits interesting behavior as a certain Riemannian metric related to the Fisher information matrix
becomes degenerate near the singularity.

\begin{figure}
\centering
    \begin{subfigure}[b]{.32\textwidth}
    \centering
        \includegraphics[scale=.5]{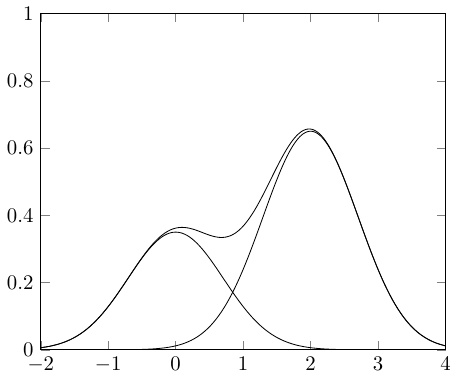}
        \caption{$a=.35,b=2$}
        \label{fig:slta}
    \end{subfigure}
    \begin{subfigure}[b]{.32\textwidth}
    \centering
        \includegraphics[scale=.5]{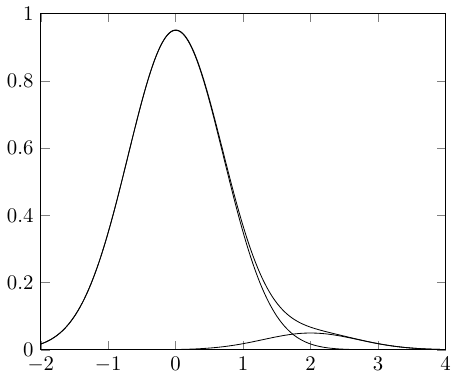}
             \caption{$a=.95,b=2$}
        \label{fig:sltb}
    \end{subfigure}
    \begin{subfigure}[b]{.32\textwidth}
    \centering
        \includegraphics[scale=.5]{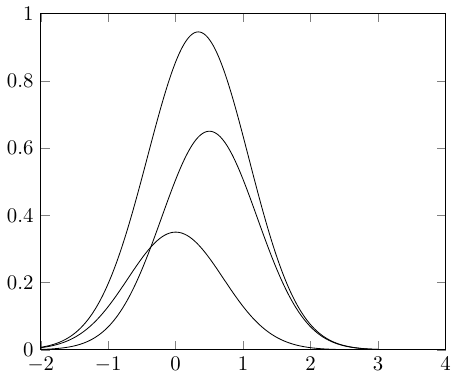}
             \caption{$a=.35,b=.5$}
        \label{fig:sltc}
    \end{subfigure}

    \caption{Simple singular statistical model with two Gaussians. The singularity arises because the second two frames correspond to different parameter values, but result in nearly equal distributions.}
    \label{fig:slt1}
\end{figure}

Consider a simple one-variable Gaussian mixture model
\begin{equation}
\label{singmodel}
p(x|a,b)=
ae^{-t^2/2}+(1-a) e^{-(t-b)^2/2},\quad
0\leq a \leq 1,
\end{equation}
consisting of a weighted sum of Gaussian distributions with standard deviation 1, and a varying mean in the second one. An illustration is shown in Figure \ref{fig:slt1}.

For any point cloud in $\mathbb{R}$, we can use the Metropolis algorithm and Bayes' rule to sample from a density proportional to
\[p(a,b|t)= \frac{p(t|a,b)p(a,b)}{p(t)},\]
starting with the uniform measure $p(a,b)=1$.
As usual, we do not need to know $p(t)$. If our point cloud is sampled from $p(t|a,b)$ for a particular choice of $(a,b)$, we would expect the resulting distribution 
to be supported near the ones we started with. This will indeed happen as expected if we chose values giving a bimodal distribution as in Figure \ref{fig:slta}. However, interesting things happen if our point cloud comes from a single Gaussian centered at the origin, because many values of $(a,b)$ correspond to that same distribution, as shown in Figures \ref{fig:sltb} and \ref{fig:sltc}.

\begin{figure}
\centering
    \begin{subfigure}[b]{.31\textwidth}
    \centering
        \includegraphics[scale=.55]{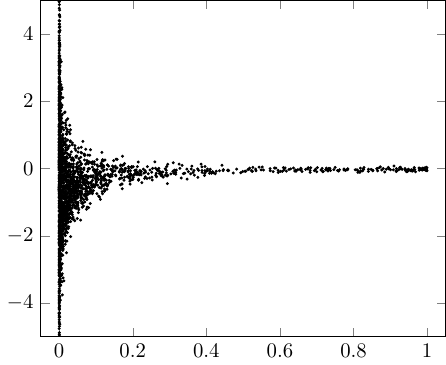}
        \caption{Many values of $a,b$}
        \label{fig:slt2a}
    \end{subfigure}
    \begin{subfigure}[b]{.31\textwidth}
    \centering
        \includegraphics[scale=.55]{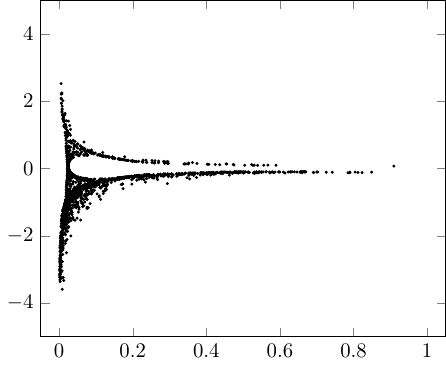}
                \caption{Samples from $\shape$}
        \label{fig:slt2b}
    \end{subfigure}
    \begin{subfigure}[b]{.33\textwidth}
    \centering
        \includegraphics[scale=.32]{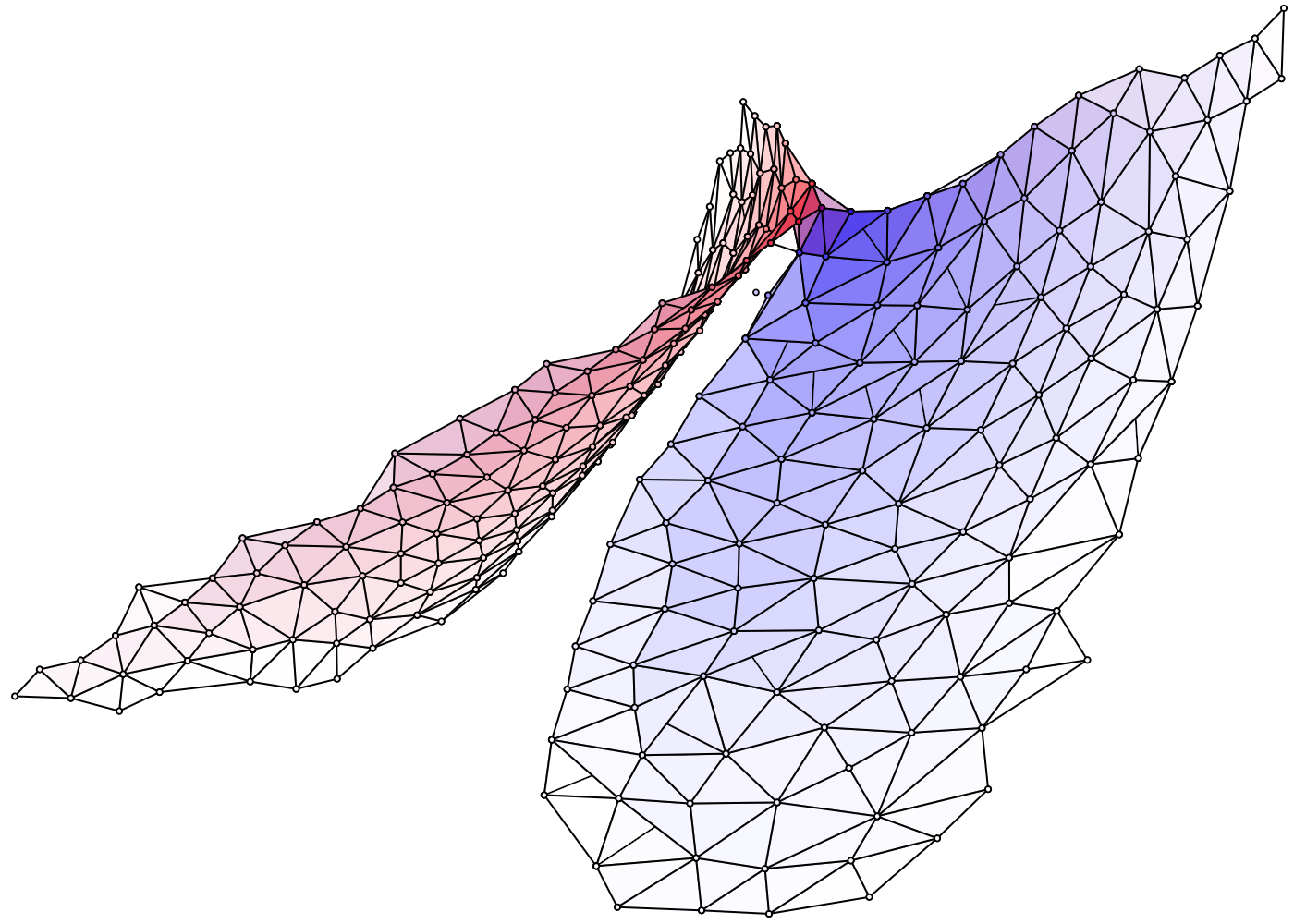}
        \caption{Alpha complex}
        \label{fig:slt2c}
    \end{subfigure}

    \caption{On the left, scatter plot of samples 
    of the model parameters $(a,b)$ using Markov chain Monte Carlo (MCMC) with respect to the distribution
    $e^{-t^2/2}$. In the middle, samples from the shape
    displayed in $\mathbb{R}^2$ using the coefficients
    of $p(x'_i)$ and the original samples
    $(a_i,b_i)$. On the right, the corresponding alpha complex $\densalphao(f,\dataset',s,d_0)$ with $s=.8$
    and $d_0=\quantile_f(.8)$.
    The red/blue colors correspond 
    to positive/negative values of $b$, respectively.
    }
    \label{fig:slt2}
\end{figure}

We generated a scatter plot of $10000$ values of the parameters $(a,b)$ associated to the singular Gaussian using 1000 training points in $\mathbb{R}$.
The results, shown in Figure \ref{fig:slt2a},
reveal the expected behavior at the singularity
at the origin $(a,b)=(0,0)$.
As in the previous section, 
building an alpha complex on the resulting
point cloud in $\mathbb{R}^2$ would be arbitrary. We no longer have a vector of predictions as in \eqref{eq:predmap}, and instead we use the vector of values of the loss function
\begin{equation}
    \label{eq:sltfisher}
(a,b)\mapsto (-\log(p(t_1|a,b)),...,-\log(p(t_{1000}|a,b))).
\end{equation}
This map has the property that the pullback Riemannian metric on $\mathbb{R}^2$
is a version of the Fisher information matrix
near the true distribution.
We then reduced the dimension down to
50 using a PCA,
in part to speed up the sampling, but more importantly
because rounding errors become a factor when evaluating
$y=p(x_i)$ at a sample. The second issue 
can easily be fixed 
by representing components of the calculation 
in log coordinates. 
We then considered the density estimator
$f:\mathbb{R}^{50}\rightarrow \mathbb{R}_+$
corresponding to the resulting point cloud
$\dataset \subset \mathbb{R}^{50}$, 
and chose a (this time arbitrary) value of 
the scale parameter of $h=.1$.

We then sampled from the shape using Algorithm \ref{alg:densland}, using the values of $s=.8$,
cutting off at $d_0=\quantile_f(.8)$.
By taking the coefficients determining
the expression 
$p(x'_i)=\sum_{j} c_{i,j} x_j$
as a convex combination of the $x_i$
and plotting the points
$\sum_{j} c_{i,j} (a_j,b_j)\in \mathbb{R}^2$ instead of
$y'_i=\sum_{j} c_{i,j} x_{j}\in \mathbb{R}^{50}$, we obtain a description of the shape in 
the original 2-dimensional plane.
Instead of corresponding to the original singularity, 
the points on the line $b=0$ were all mapped to the same point near the origin, creating the interesting pattern shown in Figure \ref{fig:slt2b}.
This is because all those values become
very close together when mapped to $\mathbb{R}^{50}$, resulting in a degenerate metric at the origin.
We then generated $\densalpha(f,\dataset',s)$ with 
$s=.8$, and plotted the projection onto a unitary subspace in $\mathbb{R}^{50}$, showing the interesting shape 
in Figure \ref{fig:slt2c}.

\subsection{Local patches in the MNIST data set}

\label{sec:patches}

In \cite{carlsson2008klein},
the authors studied the topology of a certain space of  local $3\times 3$ high intensity patches of the van Hateren data set of natural images, which was investigated earlier by Lee, Mumford, and Pederson \cite{hateren_schaaf_1998,mumford2003nonlinear}. They gave quantitative evidence using the witness complex that those patches lie along a sublocus of a parametrized
Klein bottle, called the three-circle model. 
Using a similar setup, we apply our construction to data sets of coordinate patches taken from $28\times 28$ images of handwritten digits 
from the MNIST data set \cite{deng2012mnist}.
Instead of using small $3\times 3$ patches, we project onto discrete versions of the Hermite polynomials up to quadratic order.
Using the alpha complexes,
we obtain surprisingly descriptive geometric models corresponding to different regions in the Klein bottle as one varies the digit.

A ``local image patch'' will mean an 
$l\times l$ subimage of a larger one, 
which in our case will be taken 
the MNIST data set of handwritten digits.
We will view local image patches as elements
of the vector space $V=\mat(l,l)$ using its grayscale intensity value. We have a scalar product on $V$ given by
\begin{equation}
    \label{eq:hermsp}
(A,B)=\frac{1}{2^{2(l-1)}}\sum_{i=1}^l \sum_{j=1}^l
\binom{l-1}{i-1} \binom{l-1}{j-1}
A_{i,j} B_{i,j}
\end{equation}
This inner product has a number of advantages over the usual $L^2$ product in that the weights fall off gradually near the image border. It also has the property of being nearly rotationally invariant for larger values of $l$, as the binomial coefficient approximates the Gaussian.
There is an orthonormal basis 
given by $H_{a,b}=H_a \otimes H_b$, where the $H_a \in \mathbb{R}^l$ are a discrete form of the Hermite polynomials.
They can be obtained by applying the Gram-Schmidt algorithm to the vectors of polynomial functions $v_a=(i^a)_{i=1}^l$, with respect to the one-dimensional form of \eqref{eq:hermsp}.

\begin{figure}
    \centering
  \begin{subfigure}[b]{.3\textwidth}
  \centering
\includegraphics[scale=.6,frame]{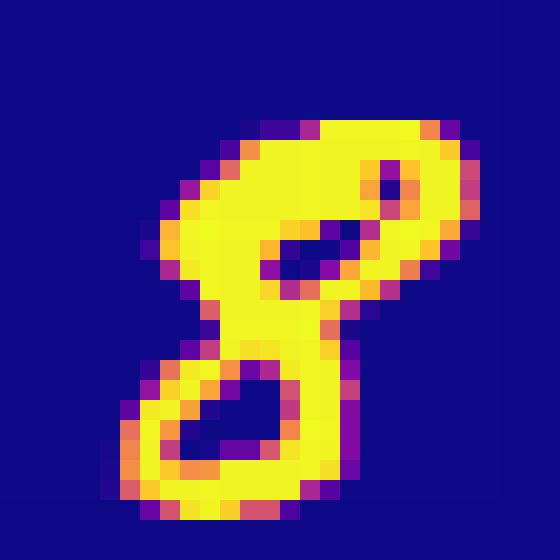}
\end{subfigure}
  \begin{subfigure}[b]{.3\textwidth}
\centering
\includegraphics[scale=1.5272,frame]{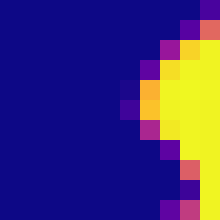}
\end{subfigure}
  \begin{subfigure}[b]{
  .3\textwidth}
  \centering
\includegraphics[scale=1.5272,frame]{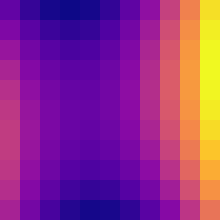}
\end{subfigure}
    \caption{A typical digit, a random $11\times 11$ image patch, and its orthogonal projection 
    onto $V_1\oplus V_2$.}
    \label{fig:hermite}
\end{figure}

We have a decomposition
\[V=V_0\oplus V_1\oplus \cdots,\quad 
V_i=\spn\{ H_{a,b}:a+b=i\}\]
as well as orthogonal projections
$\pi_i:V\rightarrow V_i$.
We will be interested in the images of image patches under the map
$\pi_{1,2}: V\rightarrow V_1\oplus V_2\cong \mathbb{R}^5$, which analogous to projecting onto low-frequency modes in Fourier analysis.
One such projection is shown in Figure \ref{fig:hermite}. For any image patch we have its norm squared $r^2=r_1^2+r_2^2$ 
where $r_i$ is the norm of its image in $V_i$. Image patches centered on the points of the 
digit would tend to have relatively high $r_2$ values, while points on the boundary would have higher $r_1$ values, due to the gradient.
The patch in Figure \ref{fig:hermite} would both have relatively high $r_1$ and $r_2$ terms.

For 50 instances of each digit, we sampled all $l\times l$ patches using the choice of $l=11$, and projected those patches onto their linear and quadratic components $V_1\oplus V_2$, to obtain a point cloud of size $50\cdot (28-l+1)^2=16200$ in $\mathbb{R}^5$. We then chose only those images whose $L^2$-norm is above a fixed number of $r\geq .3$, resulting in a subset of around 
$20\%$ of the original size.
This is analogous to the step of selecting ``high intensity patches'' 
from \cite{carlsson2008klein}. 
We then divided the remaining points by $r$ to arrive at a point cloud $\dataset_k\subset S^4$ of size a few thousand for each digit $k\in \{0,...,9\}$. 
We defined a kernel density estimator $f:\mathbb{R}^5\rightarrow
\mathbb{R}_+$ using $h=.15$,
and built $\densalphao(f,\dataset',s,d_0)$
using the same choices of $s=.5$, and $d_0=\quantile_f(.6)$ for each digit.

\begin{figure}
    \centering
  \begin{subfigure}[b]{.24\textwidth}
\includegraphics[scale=.3,frame]{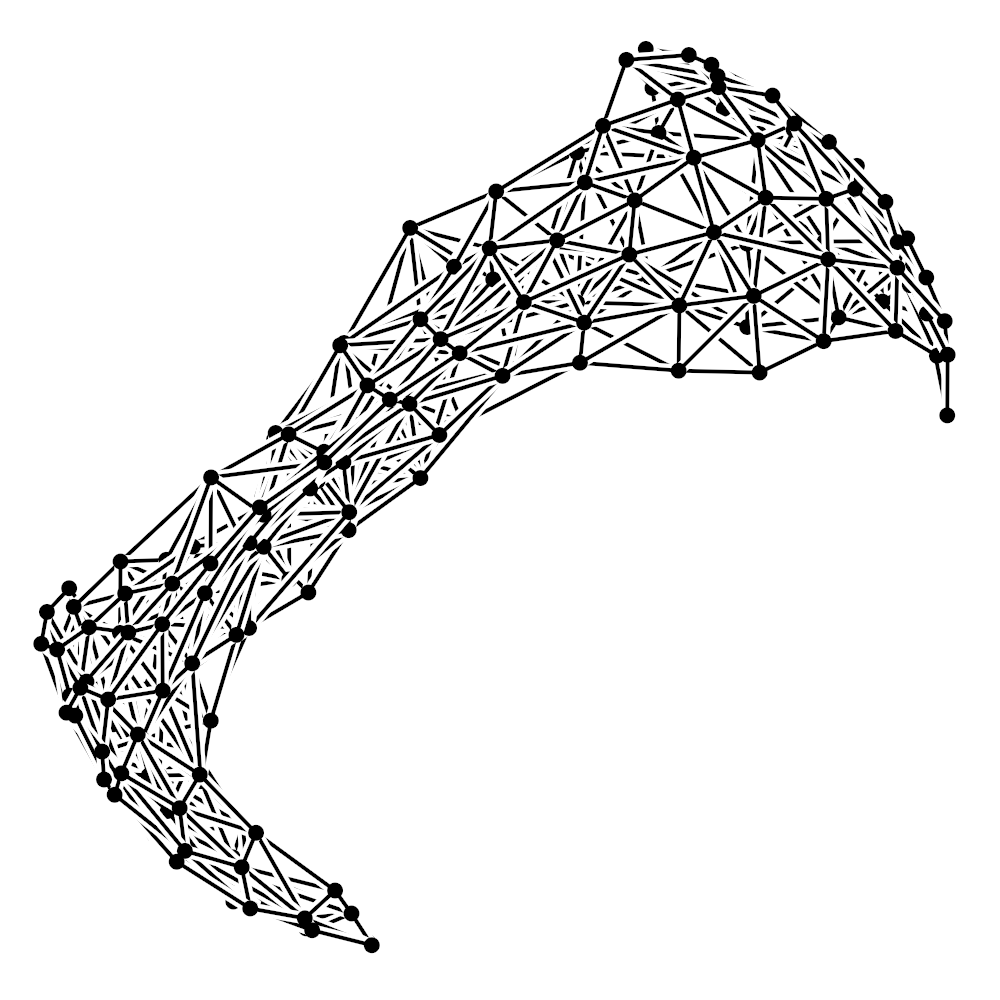}
\caption{digit=1}
\end{subfigure}
  \begin{subfigure}[b]{.24\textwidth}
\includegraphics[scale=.3,frame]{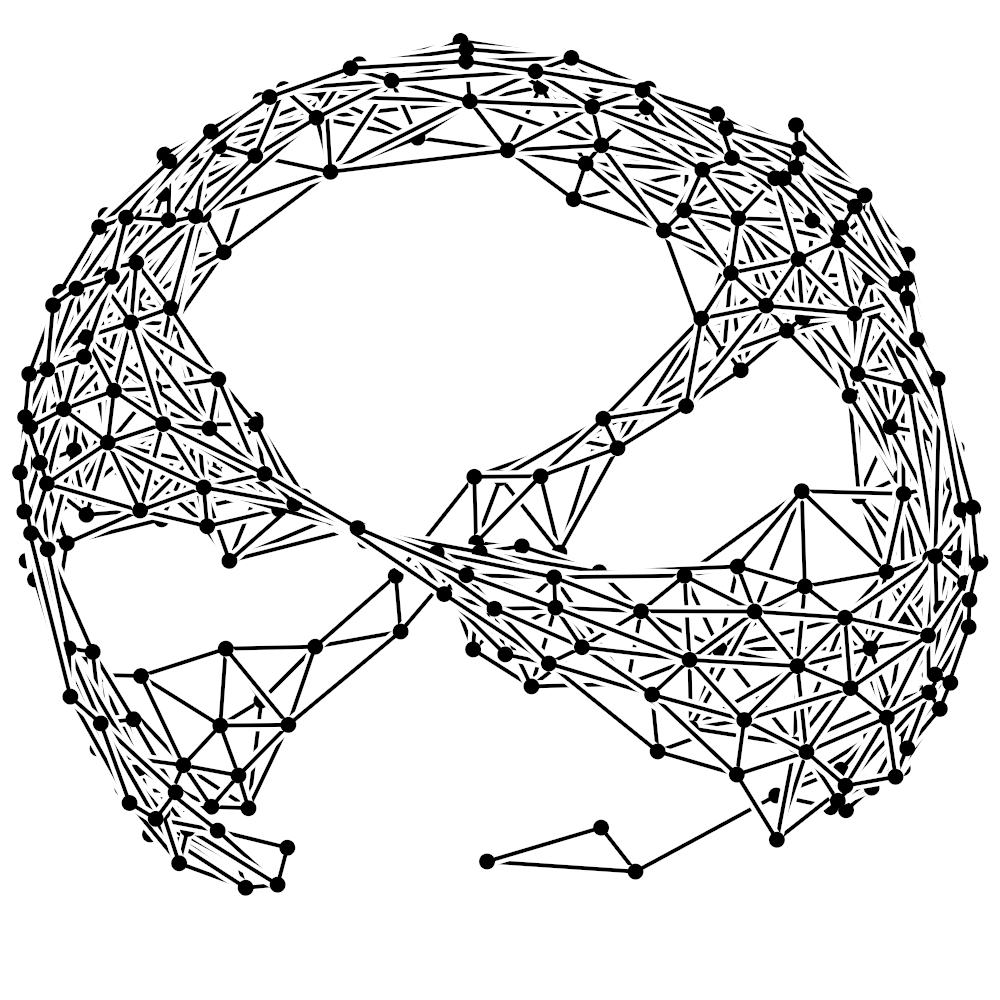}
\caption{digit=7}
\end{subfigure}
\begin{subfigure}[b]{.24\textwidth}
\includegraphics[scale=.3,frame]{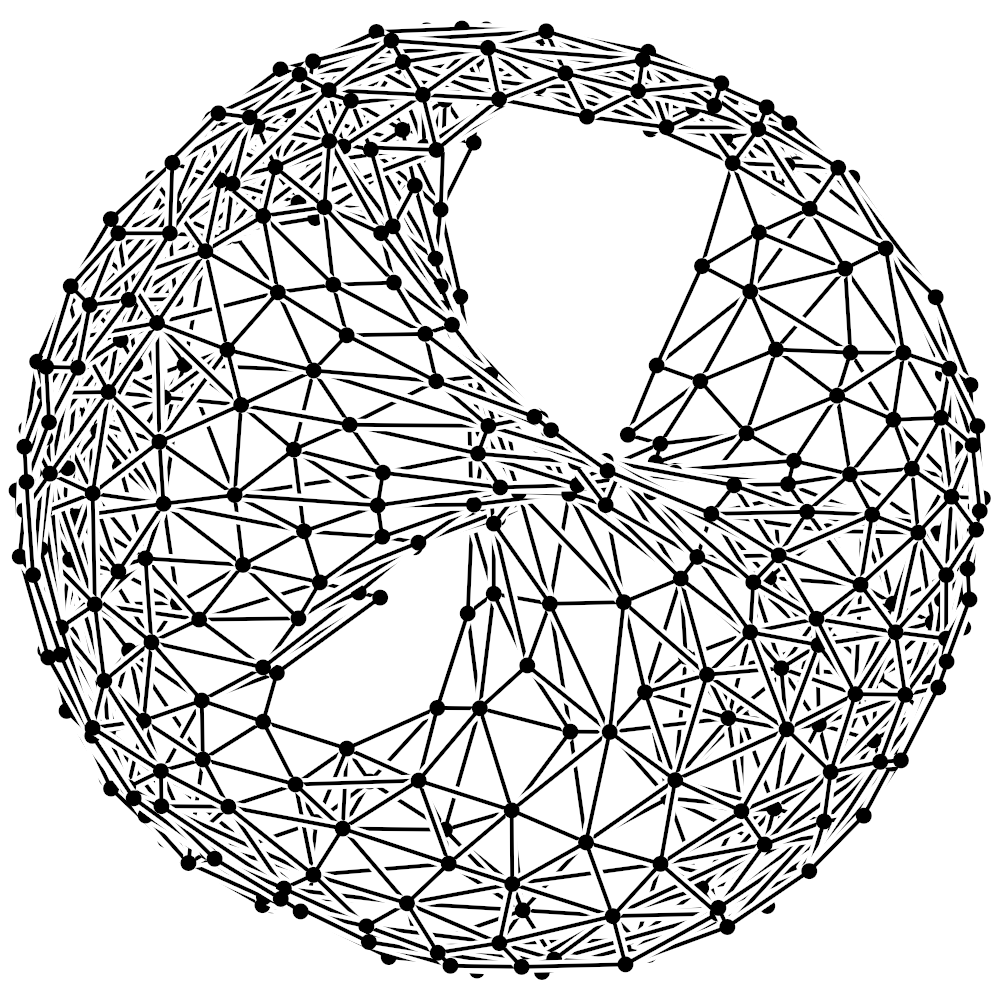}
\caption{digit=0}
\end{subfigure}
  \begin{subfigure}[b]{.24\textwidth}
\includegraphics[scale=.3,frame]{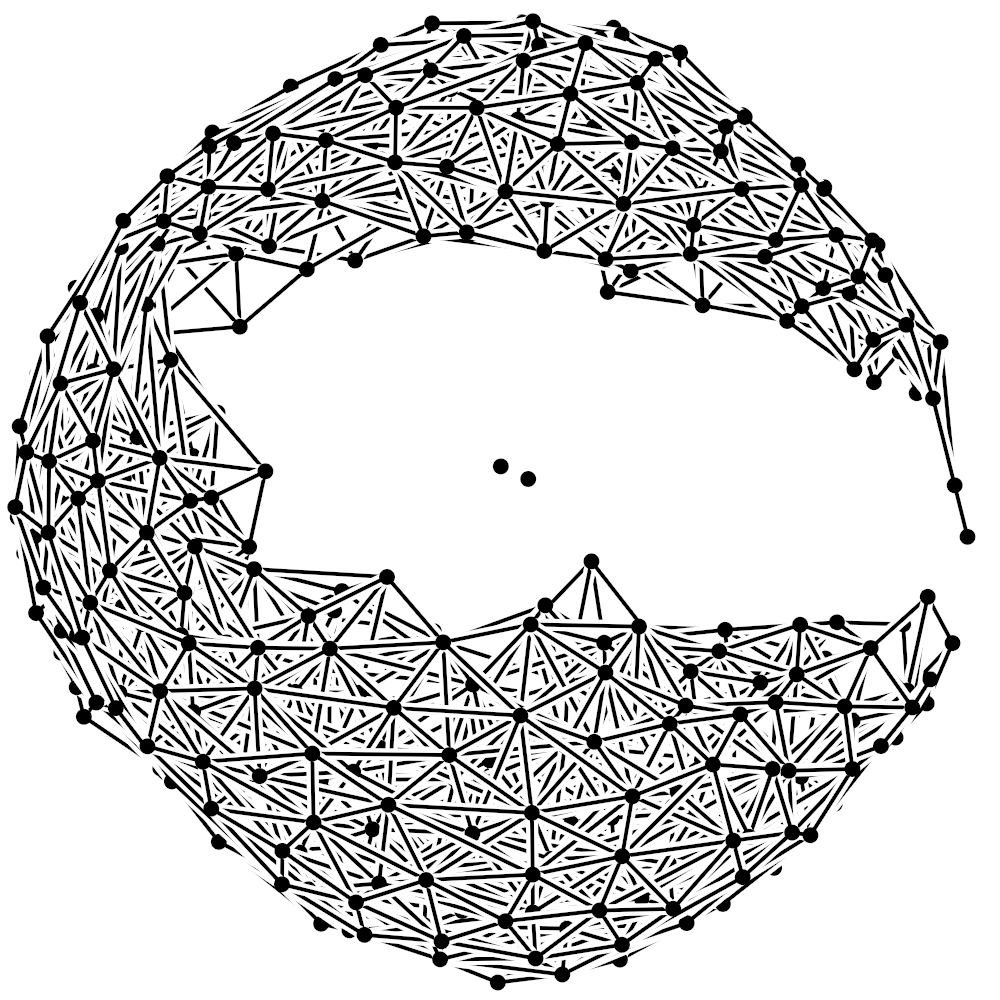}
\caption{digit=8}
\end{subfigure}
    \caption{Alpha complexes build out of the high-intensity image patches.}
    \label{fig:mnist1}
\end{figure}

The results, shown for the numbers $\{1,7,0,8\}$ in Figure \ref{fig:mnist1}, exhibited distinctive features, which can be understood in terms of the Klein bottle model.
Starting with the digit 1, we see two arcs connected by a connecting region in the center. Analyzing the images associated to each point shows that the arcs are the regions on either side of the digit, which have 
large magnitude in linear terms $r_1$. The strip lives on the digit, which has relatively high quadratic norm $r_2$. The digit 7 has a similar explanation but with two different components. These linear terms fill out the entire periphery of a circle in the the case of the digit 0, coming from both the interior and exterior. In the center of the figure, we see the high $r_2$ points
twist and connect across points on the digit,
which are only dense enough on the left and right sides because of the oval shape.
The digit 8 shows no points dominated by second order, except a small disconnected region corresponding to the two voids.

We then ran the experiment again, this time using just $r_2$ to determine intensity, and dividing by $r_2$ in place of $r$.
This has the effect of making patches with higher $r_2$ denser, thereby accentuating the second order features. This time we cut off at an intensity value of $r_2\geq .125$,
and chose $d_0=\quantile_f(.3)$, keeping the values of
$h=.15$ and $s=.5$.

\begin{figure}
    \centering
  \begin{subfigure}[b]{.24\textwidth}
\includegraphics[scale=.3,frame]{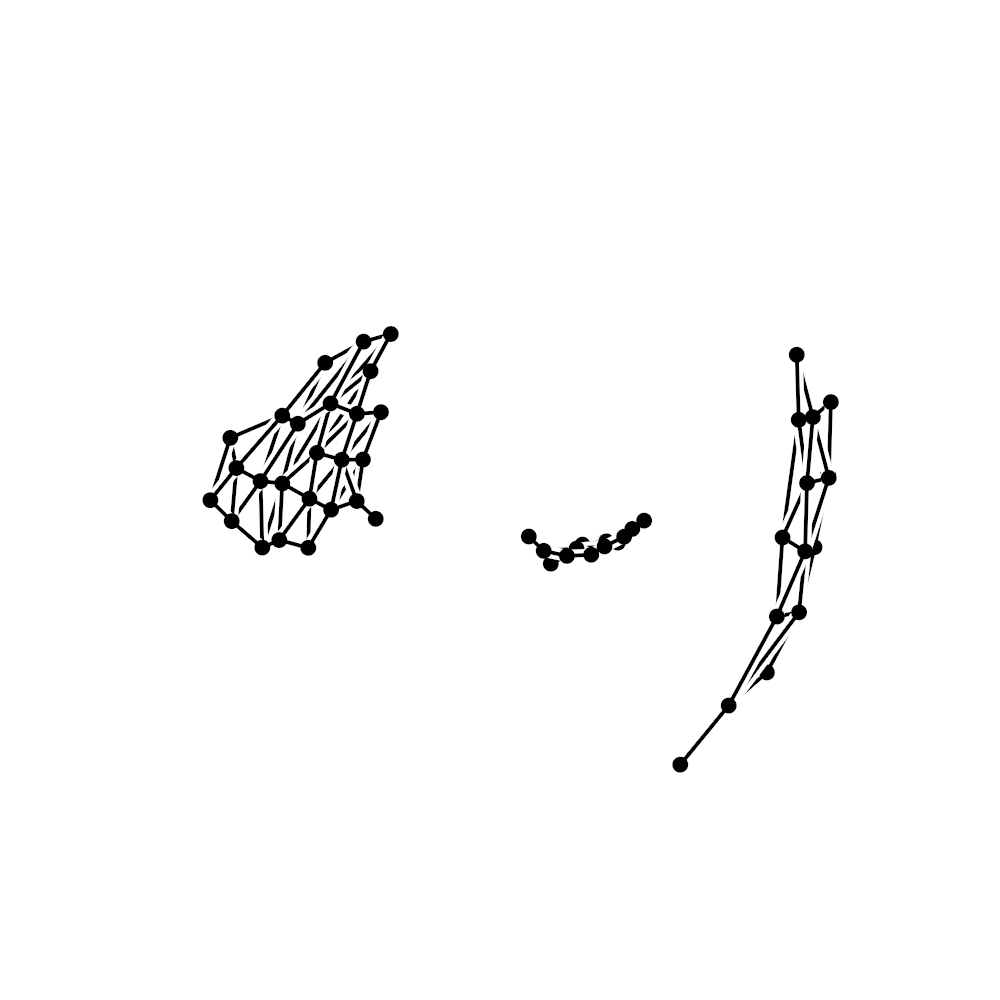}
\caption{digit=1}
\end{subfigure}
  \begin{subfigure}[b]{.24\textwidth}
\includegraphics[scale=.3,frame]{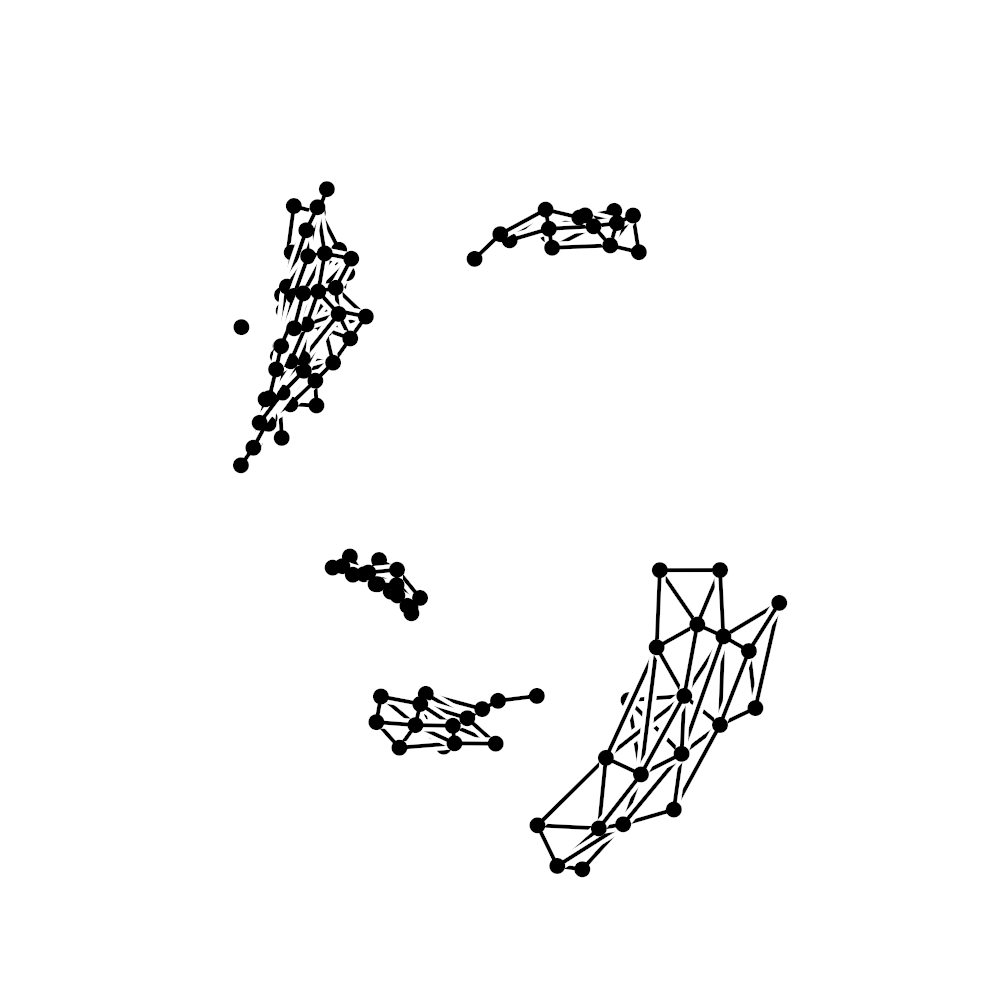}
\caption{digit=7}
\end{subfigure}
\begin{subfigure}[b]{.24\textwidth}
\includegraphics[scale=.3,frame]{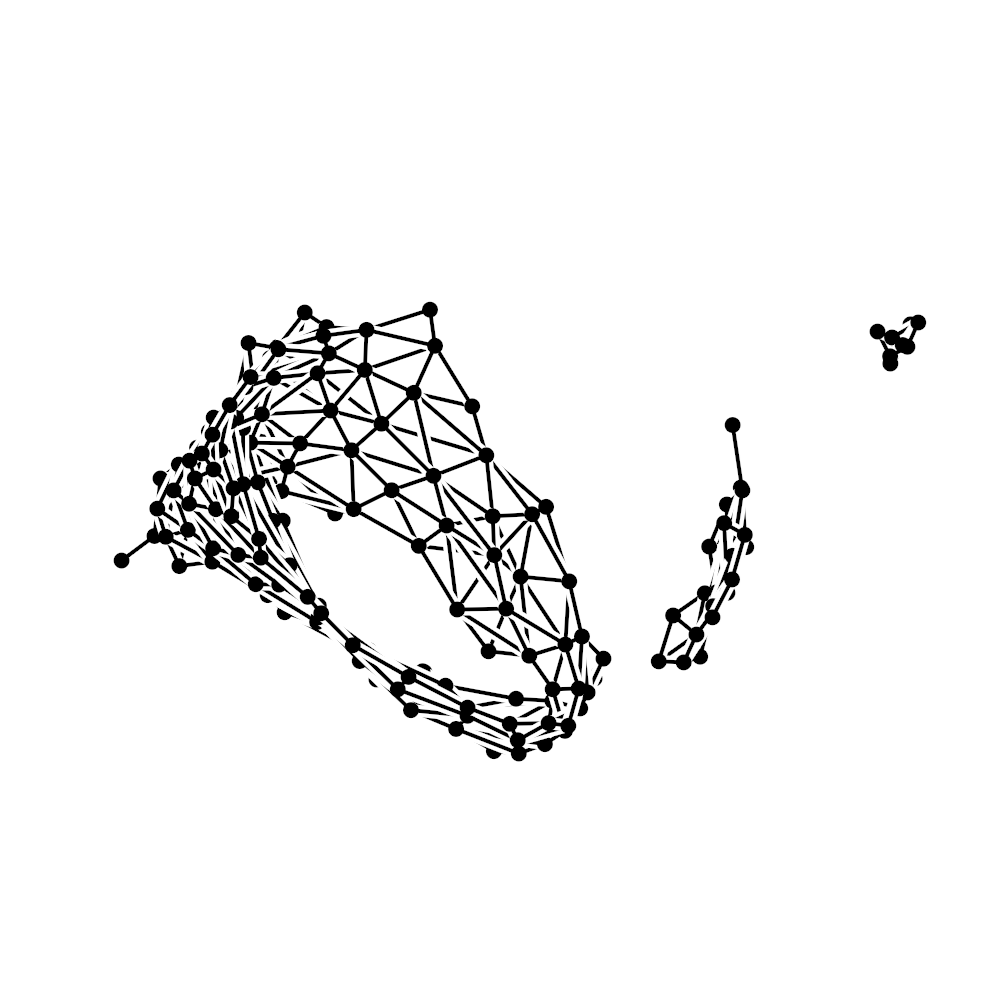}
\caption{digit=0}
\end{subfigure}
  \begin{subfigure}[b]{.24\textwidth}
\includegraphics[scale=.3,frame]{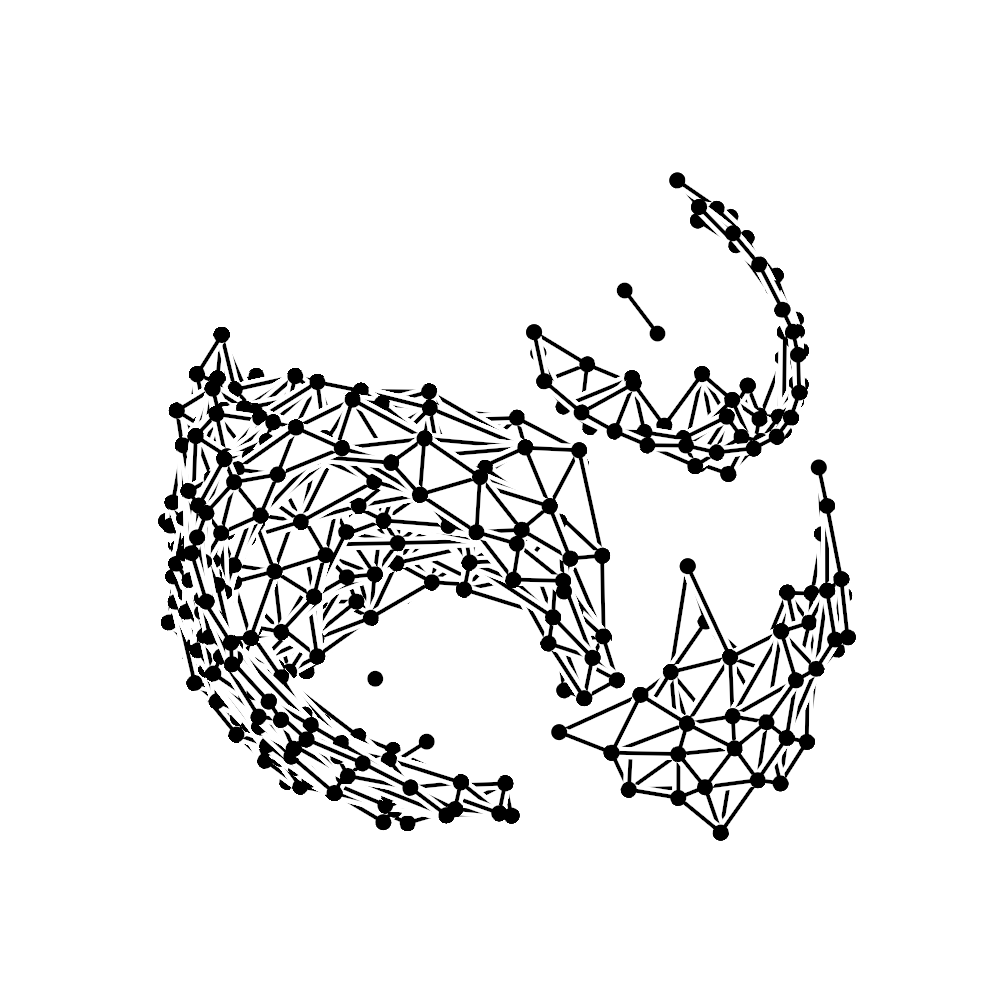}
\caption{digit=8}
\end{subfigure}
    \caption{Complexes built on the corresponding digits, normalizing only the second-order coefficients.}
    \label{fig:mnist2}
\end{figure}

The results are shown in Figure \ref{fig:mnist2}. In the digit 1, we see three connected components, corresponding to points on the digit itself, and two others which are not on the boundary, but slightly away from it on either side. For instance the patch in 
Figure \ref{fig:hermite} would be such a point.
Again, digit 7 has the same explanation, with 6 components instead of 3. In the digit zero, we have a complete M\"{o}bius strip as one traverses halfway around the digit itself. The two additional components have the same meaning as with the digit 1, appearing only on the sides because the digit is not a perfect circle. In the digit 8, we actually see 5 connected components. The main one consists of points on the digit itself, whereas the second largest ones are points on either side, as in the other digits. One of the remaining small clusters comes from the voids inside either loop, while the other represents the crossing point in the center.

\subsection{The Ising model on a graph}

\label{sec:ising}

In our final example, we consider density estimation on a simulated data set consisting of trials of the Ising model \cite{ising1925ferromagetism} on a graph with $d$ vertices, thought of as a collection of real-valued vectors in 
$\{\pm 1\}^d\subset \R^d$.
Attempting to apply kernel density estimation on the resulting point cloud directly would not yield good results, and we would not even expect to be able to correlate kernel based density at a particular state with the theoretical density determined by the energy function. We show that we can create geometric models of the density landscape as we did in Section \ref{sec:config},
by using Laplacian operator $L$
of the underlying graph to obtain
smooth versions of the spin vectors.
This suggests a way of extending the kernel density based methods to high dimensional point clouds, when the set of features comes equipped with an underlying geometry.

Let $G=(V,E)$ be a graph with $n$ vertices, 
represented by a symmetric adjacency matrix $J$, with diagonal entries being zero. In our example we will use

\hspace{1in}\includegraphics[scale=.7]{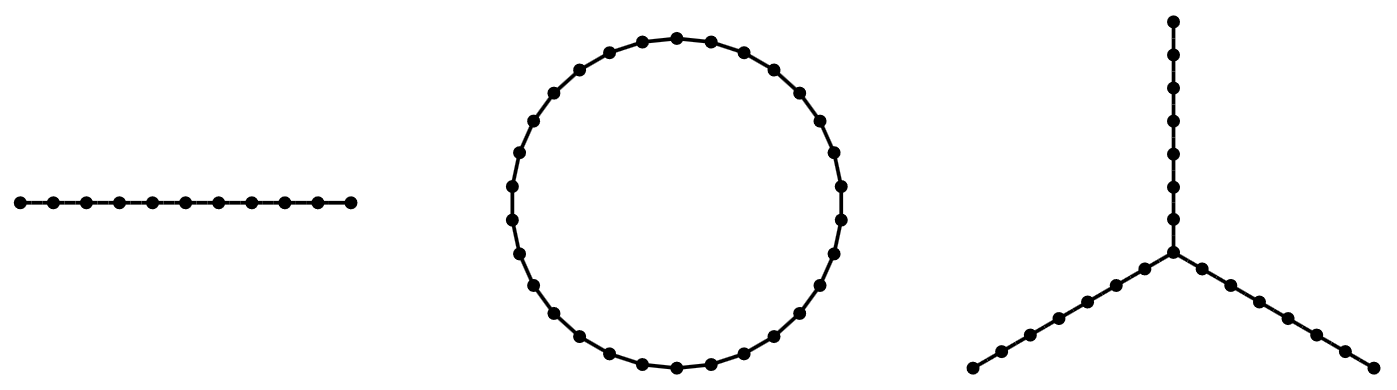}
\\
denoted $\isingint(n),\isingcirc(n),\isingflares(n)$,
where $n$ is the number of vertices.
For every discrete spin vector 
$\sigma : V \rightarrow \{1,-1\}$, 
we have the Hamiltonian energy
\begin{equation}
\label{eq:isingham}
    H_G(\sigma)=-\sum_{i,j}
J_{i,j} \sigma_i \sigma_j=
H_{\min}+2|\{(i,j)\in E: \sigma_i\neq \sigma_j\}|.
\end{equation}
Those pairs ${i,j}\in E$ for which 
$\sigma_i\neq \sigma_j$ are called transitions.
For each choice of $\beta>0$, called the temperature parameter, one seeks to sample from the
Boltzmann distribution on $\{1,-1\}^d$
given by
\begin{equation}
    \label{eq:boltzmann}
P_{\beta}(\sigma)=\frac{1}{Z_{\beta}}
e^{-\beta H(\sigma)},\quad
Z_\beta = \sum_{\sigma} e^{-\beta H(\sigma)}
\end{equation}
which is usually done using the 
single-flip Metropolis algorithm.

\begin{figure}
    \centering
    \begin{subfigure}
        [b]{.3\textwidth}
        \centering
\includegraphics[scale=.7,frame]{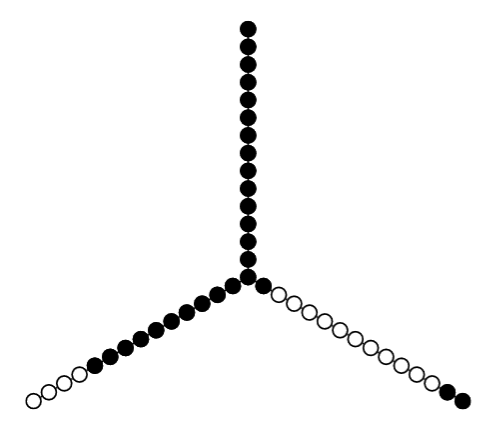}
\caption {$\sigma:V\rightarrow \{\pm 1\}$}
    \end{subfigure}
        \begin{subfigure}
        [b]{.3\textwidth}
        \centering
\includegraphics[scale=.7,frame]{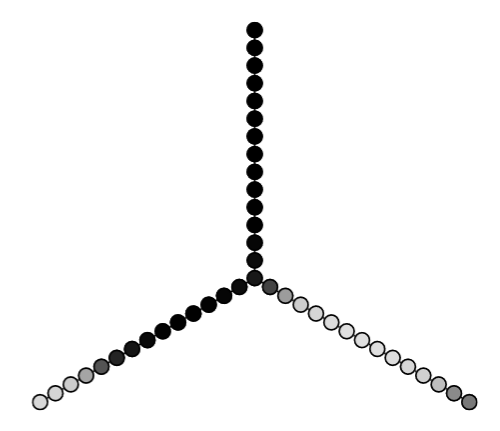}
\caption{$x=\sigma\exp(-tL)$}
    \end{subfigure}
        \begin{subfigure}
        [b]{.3\textwidth}
        \centering
\includegraphics[scale=.58]{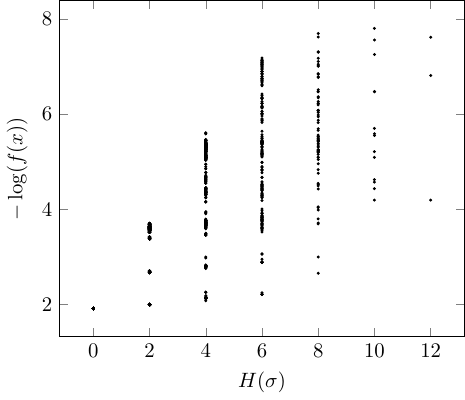}
    \end{subfigure}
    \caption{The result of blending a typical state of the Ising model simulation using the Laplacian operator. On the right, 
    a plot of the energy level versus 
    $-\log(f(x))$, when $f(x)$ is associated to the blended data set.}
    \label{fig:isingblend}
\end{figure}

For the graphs $G\in \{\isingint(30),\isingcirc(30), \isingflares(43)\}$, we simulated 
$N=20000$ states using a temperature value of 
$\beta=1.5$, and interpreted the resulting
collections of spin vectors $\{\sigma\}$
as a point clouds
$\dataset\subset \{\pm 1\}^d\subset \mathbb{R}^d$ where $d=n$ is the number of vertices. We then took a blended version 
of $\dataset$ 
using the left-normalized Laplacian operator 
$I-D^{-1}A$, where $A$ is the adjacency matrix of $G$, normalized so that the diagonal entry $A_{i,i}$ is the degree of $v_i$, and $D$ is the row-sum of $A$.
We then replaced $\dataset$
by sending
$\sigma \mapsto \sigma \exp(-tL^{t})$ with the value of $t=10$, so that the vectors are no longer $\{\pm 1\}$-valued, and considered the corresponding kernel density estimator $f:\mathbb{R}^d\rightarrow \mathbb{R}_+$ with $h=2.0$. An illustration of the result of the convolution, and the distribution of density versus 
Hamiltonian energy is shown in Figure \ref{fig:isingblend}.

\begin{figure}
\centering
\begin{subfigure}[b]{.3\textwidth}
\centering
\includegraphics[scale=.35]{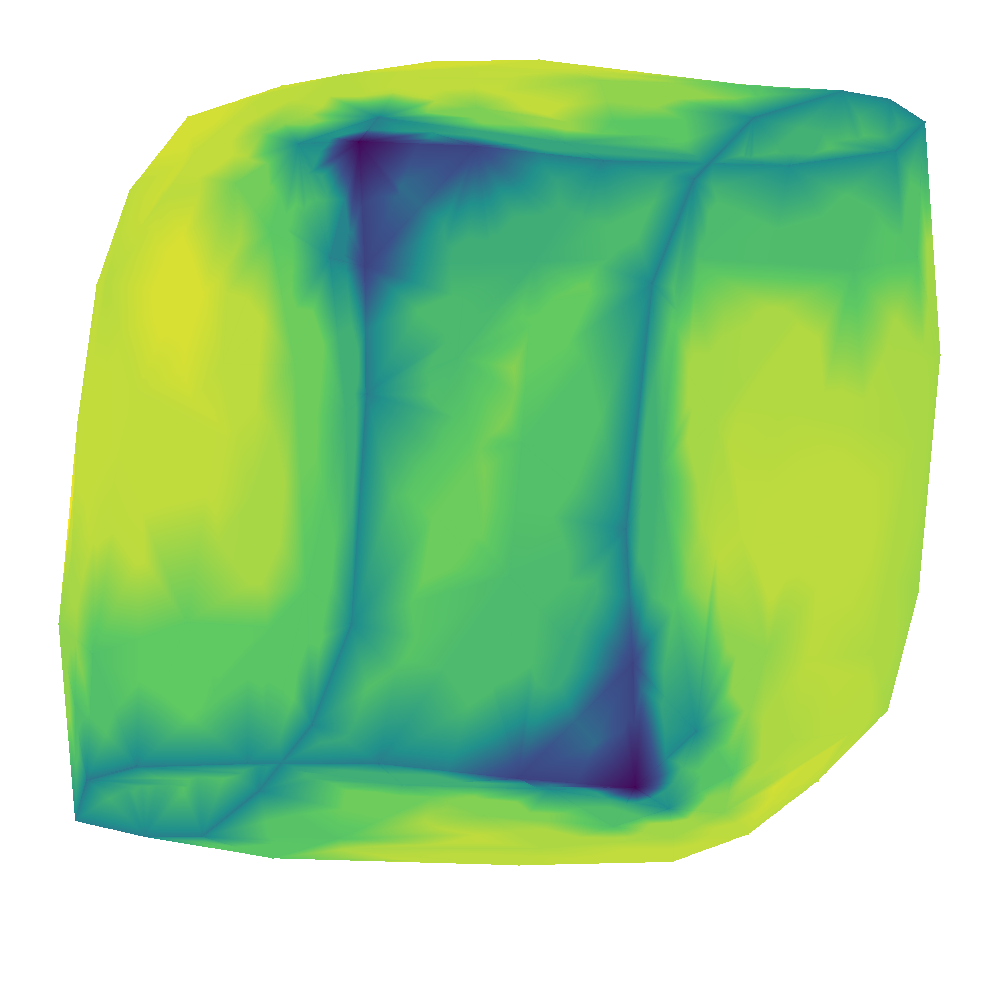}
\caption{$G=\isingint(30)$}
\label{fig:isingvisi}
\end{subfigure}
\begin{subfigure}[b]{.3\textwidth}
\centering
\includegraphics[scale=.35]{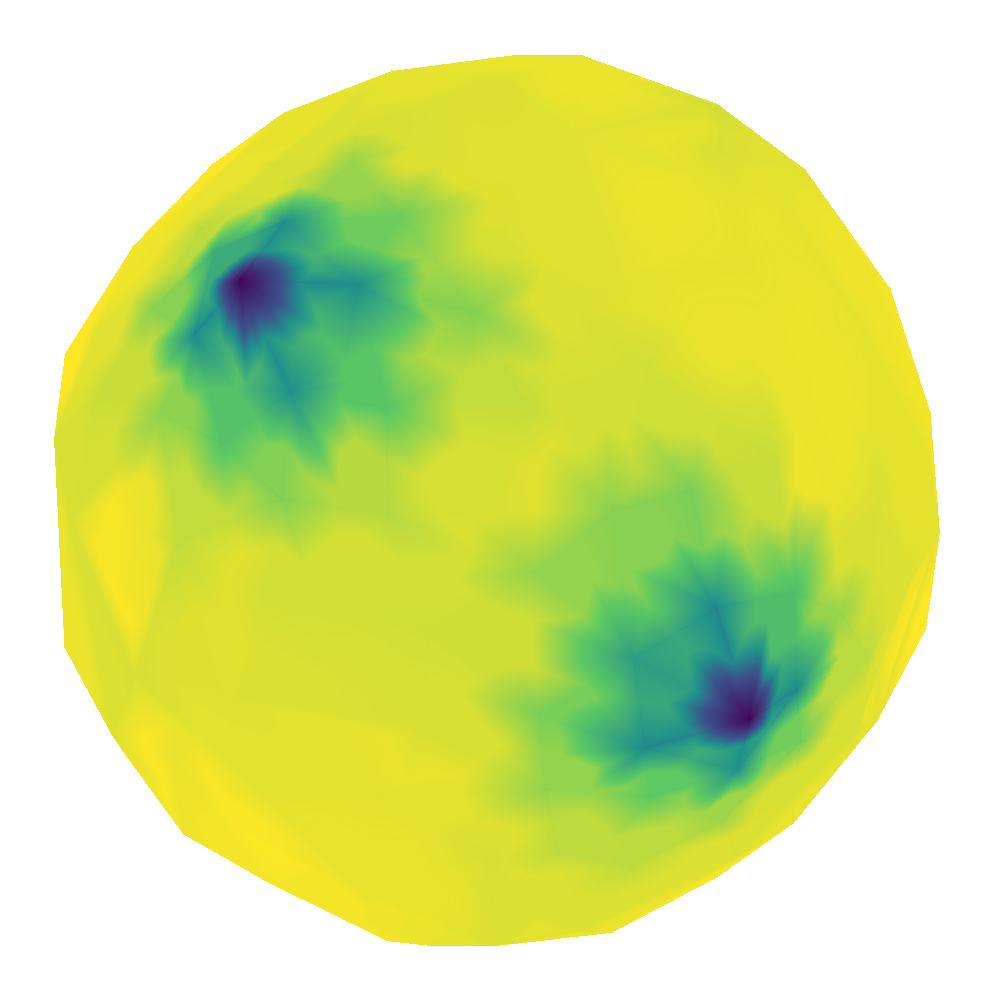}
\caption{$G=\isingcirc(30)$}
\label{fig:isingvisc}
\end{subfigure}
\begin{subfigure}[b]{.3\textwidth}
\centering
\includegraphics[scale=.35]{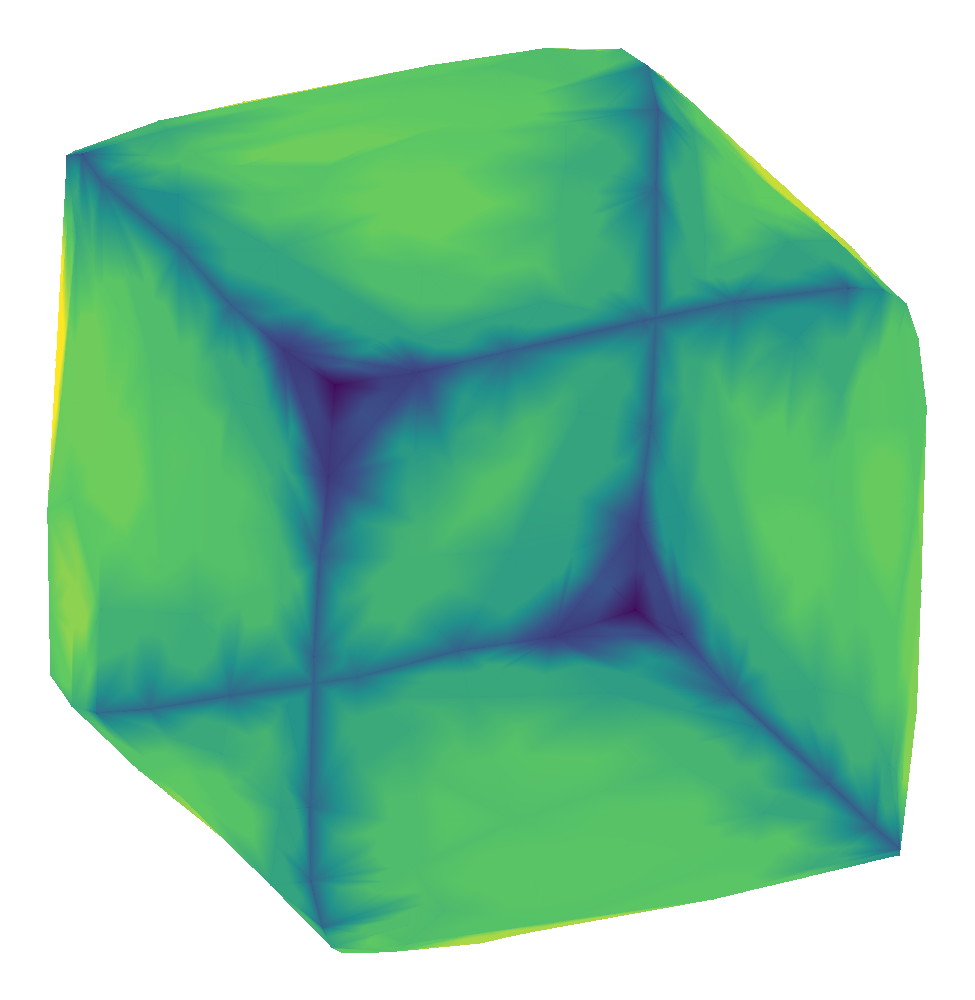}
\caption{$G=\isingflares(43)$}
\label{fig:isingvisf}
\end{subfigure}
\centering
\begin{subfigure}[b]{\textwidth}
\centering
    \includegraphics[scale=.55]{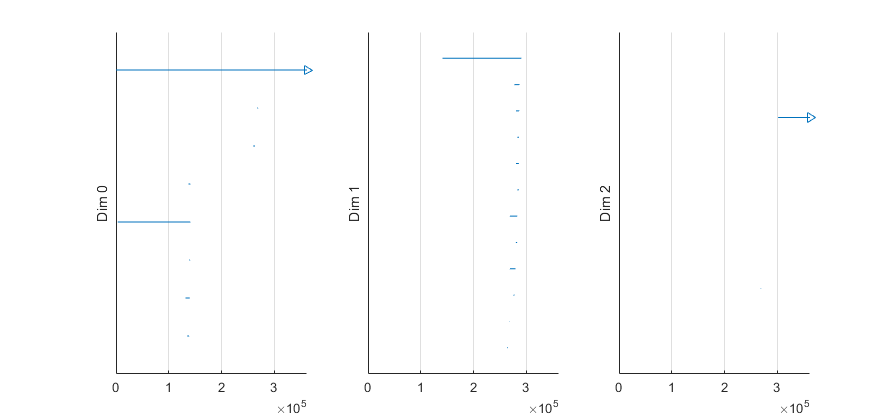}
    \caption{Barcodes associated to $G=\isingint(30)$}
\end{subfigure}
\caption{Low-dimensional projections of 
$\densalphao(f,\dataset',s,d_0)$ for the graphs
$\isingint(30)$, $\isingcirc(30)$, and $\isingflares(43)$, using $s=.5$ and $d_0=\quantile_f(.95)$. In the lower row, we have the persistence barcodes in the case of the interval.}
\label{fig:isingvis}
\end{figure}

We then computed $\densalphao(f,\dataset',s,d_0)$
for $s=.5$, and $d_0=\quantile_f(.95)$.
By taking a random 3D to 2D projection of the values in the 3 most dominant eigenvalues of $L$, we obtain a visualization of the energy landscape as in \ref{sec:config}, shown in Figure \ref{fig:isingvis}. In the first frame, associated to the interval, we have two densest types of points with energy level zero, corresponding to all spins equal to plus or minus 1, realized at the corners. We then have two curved line segments consisting of energy states with exactly one transition joining those points, starting from either side. The points shown in green/yellow correspond to states with two transitions, filling in the resulting circle to form a sphere. This is reflected in the persistence barcodes, shown in the figure.

\bibliographystyle{plain}
\bibliography{refs}

\end{document}